\documentclass[12pt, a4paper,leqno]{amsart}
\usepackage{amsmath,amsthm,amscd,amssymb,amsfonts,amsbsy}
\usepackage{latexsym}
\usepackage{txfonts}
\usepackage{exscale}

\usepackage[colorlinks,citecolor=red,pagebackref,hypertexnames=false]{hyperref}

\hypersetup{urlcolor=blue}

\calclayout
\allowdisplaybreaks

\numberwithin{equation}{section}

\usepackage{pgf}
\usepackage{color}

\newcommand{\Bk}{\color{black}}
\newcommand{\Rd}{\color{red}}

\theoremstyle{plain}
\newtheorem{theorem}[equation]{Theorem}
\newtheorem{lemma}[equation]{Lemma}
\newtheorem{corollary}[equation]{Corollary}

\theoremstyle{definition}

\theoremstyle{remark}

\newcommand{\ms}{\medskip}

\newcommand{\bP}{\mathbb{P}}
\newcommand{\bE}{\mathbb{E}}

\renewcommand{\d}{\partial}
\newcommand{\dist}{\,\mathrm{dist}\,}
\newcommand{\sm}{\setminus}

\newcommand{\wt}{\widetilde}

\newcommand{\ol}{\overline}

\newcommand{\RR}{{\mathbb{R}}}
\newcommand{\NN}{{\mathbb{N}}}
\newcommand{\ZZ}{{\mathbb{Z}}}

\newcommand{\eps}{\varepsilon}

\newcommand{\bp}{\noindent {\it Proof}.\,\,}
\newcommand{\ep}{\hfill$\Box$ \vskip 0.08in}

\newcommand{\po}{\partial\Omega}

\newcommand{\F}{\mathcal{F}}

\newcommand{\Ak}{\mathfrak{A}}

\begin{document}

\title[The landscape law]{The landscape law for the integrated density of 
states}
\author{G. David, M. Filoche, and S. Mayboroda.}
\newcommand{\Addresses}{{
  \bigskip
  \vskip 0.08in \noindent --------------------------------------
\vskip 0.10in

  \footnotesize

  G.~David, \textsc{Universit\'e Paris-Saclay, Laboratoire de Math\'{e}matiques d'Orsay, 91405, France
}\par\nopagebreak
  \textit{E-mail address}: \texttt{Guy.David@universite-paris-saclay.fr}

 \medskip
 
M.~Filoche, \textsc{Physique de la Mati\`ere Condens\'ee, Ecole Polytechnique, CNRS, Institut Polytechnique de Paris, Palaiseau, France}\par\nopagebreak
  \textit{E-mail address}: \texttt{marcel.filoche@polytechnique.edu}

\medskip

S.~Mayboroda, \textsc{School of Mathematics, University of Minnesota, 206 
Church St SE, Minneapolis, MN 55455 USA}\par\nopagebreak
  \textit{E-mail address}: \texttt{svitlana@math.umn.edu}
}}
\date{}
\maketitle

\begin{abstract} The present paper establishes non-asymptotic estimates from above and below 
on the integrated density of states of the Schr\"odinger operator $L=-\Delta+V$, using a counting function for the minima of the localization landscape, a solution to the equation $Lu=1$. 

\ms\noindent{\tt R\'esum\'e en Fran\c cais.} 
Dans cet article on \'etablit des bornes inf\'erieures et sup\'erieures sur la densit\'e d'\'etats int\'egr\'ee
pour l'op\'erateur de Schr\"odinger $L=-\Delta+V$, \`a l'aide d'une fonction comptant les minimas
de la fonction paysage, la solution de $Lu=1$ avec des conditions au bord adapt\'ees.
contexte des potentiels d\'esordonn\'es on en d\'eduit les meilleures estimations connues sur la densit\'e d'\'etats int\'egr\'ee dans le mod\`ele 
d'Anderson sur $\RR^d$. 
 \end{abstract}



\tableofcontents

\section{Introduction}\label{intro}
The density of states of the Schr\"odinger operator $-\Delta +V$ is one of the main characteristics defining the physical properties of the matter. At this point, most of the known estimates for the integrated density 
of states pertain to two {\it asymptotic} regimes, each carrying restrictions on the  underlying potentials. The first one stems from the Weyl law 
and its improved version due to the Fefferman-Phong uncertainty principle 
\cite{F}. It addresses the energies or eigenvalues $\lambda\to +\infty$ and deteriorates for the potentials oscillating at a wide range of scales. 
The second one concentrates on the asymptotics as  $\lambda$ tends to $0$ 
for disordered potentials, the so-called Lifschitz tails, and takes advantage of 
probabilistic  arguments and the random nature of the disordered potentials. The goal of 
the present paper is to establish new bounds on the integrated density of 
states via the counting function of the so-called localization landscape \cite{FM-PNAS}. The main theorem can be viewed as a new version of the uncertainty principle, which, contrary to the above, applies uniformly across the entire spectrum and covers all potentials bounded from below irrespectively of their nature.

To set the stage, let us consider the spectrum of the Schr\"odinger operator $L=-\Delta +V$ on a domain $\Omega\subset \RR^d$.  We shall assume for the time being that $\Omega$ is a cube in $\RR^d$  of sidelength $R_0\in \NN$ and make sure that the estimates that we seek 
do not depend on the size of the domain, so that we can pass to the limit 
of infinite domain whenever it is desired and appropriate. 

Assume furthermore that $V$ is a bounded non-negative function on $\Omega$ and $L=-\Delta +V$ (once again, the boundedness assumption on $V$ is, 
at this point, cosmetic: the resulting estimates do not depend on the maximum value and we can include more general potentials into consideration). We denote by $N$ the (normalized) integrated density of states of the 
operator $L$ with periodic boundary conditions on $\po$, i.e., 
\begin{equation}\label{eq2.1}
N(\mu):=\frac{1}{|\Omega|}\times \left\{\mbox{the number of eigenvalues 
$\lambda$ such that } \lambda \leq \mu\right\}.
\end{equation}
As usual, eigenvalues are counted with multiplicity.
It is known that the operator $L$ above, with periodic boundary conditions on $\po$,  has a discrete spectrum consisting of positive eigenvalues and hence, the definition is coherent.

In 1911, Hermann Weyl proposed what became later known as the Weyl law for the 
asymptotics  of $N(\mu)$, as $\mu\to+\infty$, for the Laplace-Beltrami operator with the Dirichlet boundary conditions in a bounded domain. 
In his setting, the law gives an asymptotic of a multiple of $\mu^{d/2}$ as $\mu\to +\infty$. 
Perhaps much more importantly than the result itself, it gave a general approach to the asymptotics of the density of states of an elliptic operator, and in particular,  the rule of thumb traditionally used in physics is 
\begin{equation}\label{eq1.1-bis}
N(\mu) \sim \frac{1}{(2\pi)^d |\Omega|} \iint_{|\xi|^2+V(x)<\mu} \, dxd\xi, \quad \mbox{as } \mu\to\infty.
\end{equation}
It is simultaneously impossible to list all the directions in which the Weyl law has been extended over the years and to give a sharp class of $V$ 
to which it applies, with nice control of the asymptotic errors\footnote{The estimate from above is due to Cwickel, Lieb and Rosenblum \cite{S79}.}. 
However, the oscillations of $V$ at the scales smaller than $\mu^{-1/2}$ can easily destroy the validity of the volume-counting \eqref{eq1.1-bis} for the corresponding $\mu$. In fact, the Weyl law prediction
\eqref{eq1.1-bis} fails even for systems as simple as two uncoupled harmonic oscillators, that is, the potential $V(x_1, x_2)=x_1^2+\eps x_2^2$ with a small $\eps$ (see, e.g., \cite{F}, p. 143). 

An obvious shortcoming of the ``classical" Weyl law is the emphasis
on the volume counting itself, as an eigenfunction cannot occupy an arbitrarily shaped volume in the phase space. 
This issue has been alleviated with the celebrated Uncertainty Principle of Fefferman and Phong 
ultimately reaching out to the problem of stability of matter \cite{F}. 
Instead of the volume-counting of  \eqref{eq1.1-bis}, 
Fefferman and Phong suggested to estimate the number of disjoint cubes with sidelength 
$\mu^{-1/2}$ and such that $\left(\fint_Q |V|^p\, dx\right)^{1/p} \leq C\,\mu$, smoothing the oscillations of $V$ at the correct scales. The resulting bounds on $N(\mu)$ were proved when $V$ is a polynomial
and $p=\infty$ in \cite{F} and for $V\geq 0$ in a suitable reverse H\"older case  by Shen \cite{Shen-TAMS, Shen-Duke}, and were also extended to 
estimates on a number of negative eigenvalues for general $V\leq 0$.  Overall, these ideas have brought a number of fascinating results -- their goals and achievements, stemming from a new diagonalization of pseudodifferential operators, are beyond the scope of our review. But in the particular context of interest in this paper, they also 
\Rd
fall
\Bk 
short in some respects. First,  searching for the aforementioned collection of optimal cubes for every $\mu$ can be computationally very challenging. Secondly, and this is exactly the reason for the restrictions on the potential and/or asymptotic nature of the results, the sharp estimates from above and below for positive potentials are only available when $V$ behaves not too violently at the corresponding scales. This rends them formally inapplicable for the Anderson or other disordered potentials, and more generally whenever $V$ is very different from its average on a cube. The Landscape Law proposed in this paper addresses both of these issues. The landscape ``determines" the correct cubes and exhibits precisely the correct oscillation, in some sense creating a perfect effective potential  
for the Fefferman-Phong-type counting from any initial $V$.  

above which would be desirably close to the estimate from below is challenging and requires different techniques.

landscape, which yielded astonishingly precise non-asymptotic estimates on the density of states for both periodic and certain Anderson-type potentials throughout multiple numerical and physical experiments  \cite{FM-PNAS, ADFJM-SIAM, ADFJM-PRL}. However, so far no rigorous mathematical results have supported these findings and, in particular, it was not clear what are the exact bounds, what is the range of potentials  to which the theory could be applied, whether the results are generic or governed by the 
particular choice of examples, whether one can truly furnish localization 
landscape theory in the context of Anderson localization.  In the present 
paper we prove that a  counting function arising from the landscape provides sharp estimates from above and below on the density of states {\it for any non-negative potential} in the Schr\"odinger operator.  As a by-product,  we derive new estimates on the integrated density of states for the Anderson-type potentials.  However, the latter is only a particular instance of our theory -- our  main results are deterministic.

The concept of {\it localization landscape} was pioneered by the second and third authors of the present paper in \cite{FM-PNAS}. 
The landscape is the
solution to $(-\Delta+V)\,u=1$, with the same boundary conditions as the original operator in question. When applied to the Laplacian rather than the Schr\"odinger operator and equipped with the Dirichlet boundary conditions, the landscape is  nothing else than  the classical torsion function, however, its role in our theory and its character in the presence of a potential are very different, and we will continue using the {\it landscape} terminology which seems to be more illustrative under the circumstances.

First numerical \cite{ADFJM-PRL} and then rigorous mathematical results \cite{ADFJM-CPDE} have demonstrated the relationship between the landscape 
and the location and shape of localized eigenfunctions, including the pattern of their exponential decay. 
One of the key observations underpinning these works is that the operator 
$L=-\Delta +V$ has exactly the same spectrum as a conjugated operator 
$$-\frac {1}{u^2}\, {\rm div} \, u^2 \,\nabla +\frac 1u$$ 
which brings up $1/u$ as an effective potential.  This is a 
consequence of the identity
\begin{equation}\label{eq1.2}
\int  |\nabla f|^2 + Vf^2\,   dx  
= \int  u^2 \left| \nabla \,\Bigl(\frac fu\Bigr)\right|^2 +\frac1u\, f^2  \, dx,
\end{equation}
valid for all $f$ in the corresponding Sobolev space $W^{1,2}(\Omega)$ and proved in \cite{ADFJM-CPDE}. However, not only $1/u$ plays the role of a potential, but it exhibits decisively better properties than the original $V$. The reduced kinetic energy, which is the first term on the right-hand side of of \eqref{eq1.2}, is small in many typical examples, at least at the bottom of the spectrum, and hence $1/u$ ``absorbs" the information about both kinetic and potential energy of the original system, in some sense, yielding a stronger form of the Uncertainty Principle than those discussed above.

\begin{figure}[htbp]
\includegraphics[width=3in]{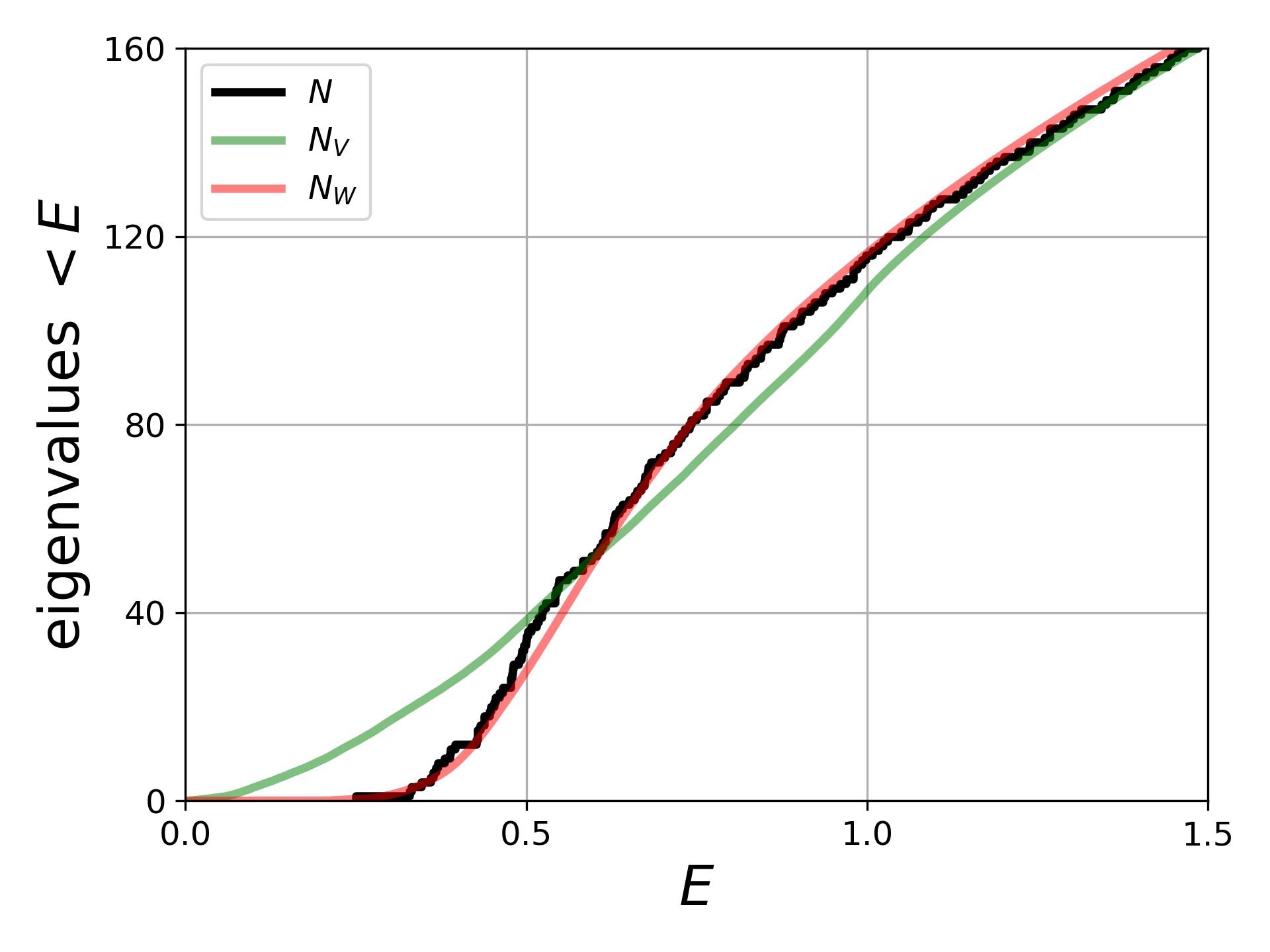}
\caption{\cite{ADFJM-SIAM} The IDOS $N$ (in black), the original Weyl law 
approximation 
$N_V$ from the right-hand side of \eqref{eq1.1-bis} (in green), 
and the approximation using the landscape function, $N_W$, $W=1/u$, from the right-hand side of \eqref{eq1.4} (in red)  for a random  uniform potential in one dimension on an interval of length 512. 
The quantities are not normalized by volume.}
\label{fg:weyl1d} 
\end{figure}

Motivated by these considerations, we were led to investigate the information about the spectrum 
of $L$ encoded in $1/u$, and the numerical experiments brought  surprising results, in fact,  exceeding original expectations \cite{ADFJM-PRL, ADFJM-SIAM}. In generic samples of Anderson-type potentials in finite one- and two-dimensional domains one could observe two strongly emerging patterns. First, the eigenvalues at the bottom of the spectrum are essentially dimensional multiples of local minima of $1/u$.  That is, independently of the potential, we observe an almost equality
\begin{equation}\label{eq1.3}\left(1+\frac d4\right) \left(\min \frac 1u\right)_j \sim \lambda_j
\end{equation}
where the eigenvalues and minima are indexed in nondecreasing order.
Secondly, a version of the Weyl law governed by the potential $1/u$ 
\begin{equation}\label{eq1.4} 
N(\mu) \sim \frac{1}{(2\pi)^d \,|\Omega|} \iint_{|\xi|^2+\frac{1}{u(x)}<\mu} \, dxd\xi
\end{equation}
yields, contrary to \eqref{eq1.1-bis}, an approximation of the density of 
states throughout the spectrum, for all values of $\mu$, albeit working a 
little worse than minima \eqref{eq1.3} at the very bottom. 
Figure~\ref{fg:weyl1d},  taken from \cite{ADFJM-SIAM}, shows the advantage of using the landscape rather than the original $V$ in the predictor \eqref{eq1.4}. 

Both observations have been immediately adopted by physicists,  for Schr\"odinger and Poisson-Schr\"odinger (Hartree-Fock) systems~\cite{LED1, LED2, Schulz-PRB, JApplPhys}, and for Dirac equation~\cite{Dirac}; however, even rigorous mathematical conjectures remained beyond reach, particularly if aiming for non-asymptotic statements. Indeed, one can rather easily construct counterexamples about taking  \eqref{eq1.3} or \eqref{eq1.4} as near identities~\cite{ComtetPRL, ReplyPRL}, and the numerical evidence  was initially restricted to dimensions 1 and 2, either Anderson-type 
or periodic potentials, and reasonably small domains, especially in dimension 2. The latter point, in particular, could raise doubts on the applicability of these approximations in the limit of infinite domain.

The present paper is the first mathematical treatment of a rigorous connection between the landscape function and the eigenvalues of $L$ in the entire range of $\lambda$. 
We show that a counting function of the minima of $1/u$ yields sharp deterministic estimates {\it from above and below} on the integrated density of states, without restrictions on the underlying potential. 

Passing to the statements of the results, recall that $\Omega$ is a cube in $\RR^d$ of sidelength $R_0\in \NN$. 
For any $r>0$ such that $R_0$ is an integer 
multiple of $r$, 
we denote by $\{Q\}_r$ a disjoint collection of cubes of sidelength $r$, such that every 
$Q_r$ is contained in $\Omega$ and $\bigcup_{Q\in \{Q\}_r} \bar Q = \bar\Omega$. 
Our cubes are always open unless stated otherwise. We shall work with functions satisfying periodic boundary conditions on $\partial\Omega$ and, slightly abusing the notation, will often identify $\Omega$ with 
the torus $(\RR/R_0 \,\ZZ)^d$.
As in the beginning of the introduction, $V$ is a bounded nonegative function on $\Omega$, $L=-\Delta +V$ is the Schr\"odinger operator on $\Omega$, which we take with the periodic boundary conditions, and the integrated density of states is defined by \eqref{eq2.1}. Going further, let $u$ 
be the solution to $Lu=1$ on $\Omega$, also with periodic boundary conditions. Then it is known (and easy to prove) that $u$ is positive and bounded, and we define 

\begin{equation}\label{eq2.2}
N_u(\mu):=\frac{1}{|\Omega|}\times \left\{\mbox{the number of cubes } Q\in \{Q\}_{\kappa\,\mu^{-1/2}} \mbox{ such that } \min_Q \frac 1u \leq \mu\right\},
\end{equation}
where by convention $1\leq \kappa<2$ (depending on $\mu$) is the smallest 
number such that  $R_0$ is an integer multiple of $\kappa \mu^{-1/2}$. 

\begin{theorem}[The Landscape law]
\label{t2.1} Retain the definitions above. There exist constants 
$C_i$, $i=1,...,4,$ depending on the dimension only, such that 
\begin{equation}\label{eq2.4}
C_1  \alpha^{d} N_u (C_2 \alpha^{d+2} \mu) - C_3N_u (C_2 \alpha^{d+4} \mu) \leq N(\mu) 
\leq N_u (C_4 \mu) 
\end{equation}
for every $\alpha<2^{-4}$ and every $\mu>0$. 
\end{theorem}

The strength of Theorem~\ref{t2.1} lies in its generality compared to all 
previously available results: 
\begin{itemize}
\item Theorem~\ref{t2.1} is not asymptotic, the estimate \eqref{eq2.4} is 
valid 
throughout the spectrum, with constants independent of $\mu$.
\item The constants in \eqref{eq2.4} do not depend on smoothness or oscillations of $V$, nor on the 
possible probability law beyond its construction (or lack of thereof), nor, in fact,  on the $L^\infty$ norm of $V$ 
or the size of the domain $R_0$. 
If one allows the dependence on $\|V\|_{L^\infty(\Omega)}$, the situation 
for large $\mu$ is 
of course rather trivial (both the density of states and $N_u(\mu)$ roughly behave as those of the Laplacian), and similarly the scales bigger than $R_0$ would be easy to handle. 
We emphasize the lack of dependence on any of these parameters, 
which makes it possible to apply the theorem 
to the limit of an infinite potential or an infinite domain.
\end{itemize}

Looking at \eqref{eq2.4}, one obviously faces the question of the polynomial correction
in the estimate from below. And indeed, in applications \eqref{eq2.4} 
often transforms into the
even stronger estimate 
$$N_u (C'_2 \,\mu) \leq N(\mu) \leq N_u (C_4\, \mu)$$
by taking $\alpha$ small. There are (at least) two mechanisms to achieve this, which are fortunately roughly complementary. The first one  is to prove a doubling condition for the landscape $u$. 

\begin{theorem}[The doubling case]
\label{t2.1-bis} Retain the definitions above. If, in addition, $u^2$ is a doubling weight at relatively small scales, specifically,  if there is a constant $C_D \geq 1$ such that 
\begin{equation}\label{eq2.5}
\int_{Q_{2s}} u^2\, dx\leq C_D \left( \int_{Q_{s}} u^2\, dx  +s^{d+4}\right) 
\end{equation}
for every cube $Q_s$ of sidelength $s>0$ 
then
\begin{equation}\label{eq2.6}
N_u (C'_2 \,\mu) \leq N(\mu) \leq N_u (C_4\, \mu)
\ \hbox{ for every $\mu>0$,}
\end{equation}
where $C_4$ is as in Theorem~\ref{t2.1} and $C'_2$ depends only on $C_D$ and the dimension.
\end{theorem}

In the doubling condition and everywhere below, 
we interpret $u$ as a function on the torus, that is, if the cubes intersect the boundary, 
it is understood that one uses the periodic extension of $u$. 

There is a  certain dichotomy between  the range of applicability of Theorem~\ref{t2.1-bis} 
and its complement, in particular, disordered systems. 
Notice that  \eqref{eq2.4}  transforms into \eqref{eq2.6} if $N_u(\mu)$ decays sufficiently fast as $\mu$ tends to $0$. 
This would not be the case, e.g., in the realm of periodic potentials, when one expects that both the integrated 
density of states and $N_u (\mu)$ behave as $\mu^{d/2}$. Fortunately, in this case $u^2$ is a doubling weight, \eqref{eq2.5} is satisfied, and hence 
we can directly apply Theorem~\ref{t2.1-bis}. 

A similar situation 
\Rd
occurs
\Bk 
when $V$ is sufficiently well-behaved. For instance, for $d\geq 3$, 
if $V$ satisfies the Kato condition
\begin{equation}\label{Kato}
\sup_{z\in \RR^d, \, r>0} \int_{B_r(z) \cap \Omega} \frac{V(x)}{|x-z|^{d-2}} \, dx =:K <\infty,
\end{equation}
then \eqref{eq2.5} is verified and hence, the integrated density of states satisfies \eqref{eq2.6} directly. 
This can be seen as a combination of results from Theorem~1.3 in \cite{Kurata},
which guarantee that for non-negative supersolutions to $(-\Delta +V) u \geq 0$ there exists $\delta>0$ such that $u^\delta$ is doubling, and classical Moser inequalities for subsolutions to $-\Delta u\leq 1$, 
which allow one to bound $\sup_{Q_s} u$ by $\left( \fint_{Q_{s}} u^2\, dx\right)^{1/2}+r^2$ (cf. \cite{HL}, 
Theorem 4.14). We observe that this includes, on finite domains, even singular potentials weaker than $1/|x|^2$, but as usual, one has to pay attention to emerging constants: if \eqref{Kato} is used, the resulting constant in \eqref{eq2.6} \Rd
may
\Bk 
depend on $K$, which might or might not be suitable for the problem at hand. 
In fact, if $V$ is regular itself,  \eqref{eq2.5} could be easier to check directly, without involving \eqref{Kato}, 
but for now let us move to the case when \eqref{eq2.5} 
\Rd
can fail. 
\Bk 

could display  pure point spectrum and exponentially decaying eigenfunctions. A certain pre-runner of Anderson localization (in fact, a simpler phenomenon of rare big regions) manifests itself through the so-called Lifschitz or Urbach tails: as $\mu\to 0$,  $N(\mu)$ behaves roughly as $e^{-c 
\mu^{-d/2}}$ contrary to the more usual behavior $\mu^{d/2}$ observed in non-disordered systems (compare to the Weyl law above). We underline that 
this, once again, is an asymptotic result, now at the edge $\mu \to 0$, with a limited understanding of errors and the range where the asymptotic is precise.

A typical example of  potential that destroys \eqref{eq2.5} is any of the 
Anderson-type  potentials. The latter is a subclass of disordered potentials 
where  $V$ is, for instance, a linear combination of bumps with random amplitudes taking values between 0 and 1 according to some probability law. It is a setting of the Anderson localization -- a famous phenomenon when such a 
system, in the limit of an infinite domain, could display  pure point spectrum and exponentially decaying eigenfunctions. 
We shall see that in this case, although \eqref{eq2.5} fails, fortunately 
$N_u(\mu)$ has exponential growth
as $\mu \to 0$, and hence
\eqref{eq2.4} implies \eqref{eq2.6} because the exponential behavior suppresses polynomial corrections. In the terminology of \cite{PasturBook}, such is the situation near fluctuation boundaries generally exhibited by Schr\"odinger operators with random (disordered) potentials.  Hence, any fluctuating boundary would yield \eqref{eq2.6}. Here we just isolate one result.

\begin{theorem}\label{t1.11}
Retain the definitions of $\Omega$ and $L$ from Section~\ref{intro}.  

Let $\varphi\in C_0^\infty(B_{1/10}(0))$
be a nontrivial bump function supported in the
ball centered at 0 of radius $1/10$, with $0\leq \varphi\leq 1$,  and set 

$$
V=V_{\omega} (x)=\sum_{j\in \ZZ^d} \omega_j \varphi(x-j) \quad \hbox{for } x\in \Omega,
$$
where the $\omega_j$ are i.i.d. variables taking values in $[0,1]$, 
whose probability distribution
$$F(\delta)=\bP\{\omega \leq \delta\}, \quad 0\leq \delta\leq 1,$$
is not trivial, i.e., not concentrated at one point, and such that $0$ is 
the infimum of its support. 
Denote by $N^E_u(\mu)=\bE\, N_u(\mu)$ the expectation of the counting function of minima 
of $1/u$, as defined in \eqref{eq2.2} and by $N^E(\mu)=\bE\, N(\mu)$ the expectation of the 
density of states, as defined in \eqref{eq2.1}. 

Then there exist constants $C_5, C_6 >0$ depending
on the dimension and the expectation of the random variables $\omega_j$ %
added an s and the name
only, and a constant $C_4>0$, depending on the dimension only, such that 
\begin{equation}\label{eq3.33-bis}
C_5N_u^E (C_6\, \mu) \leq N^E(\mu) \leq N_u^E (C_4\, \mu),
\end{equation}
for every $\mu>0$.
\end{theorem}

Since $0$ is the infimum of the support of $F$, we have $F(\delta) > 0$ for $\delta >0$;
also, the measure is not a Dirac mass at the origin, so $\lim_{\delta \to 
0} F(\delta) < 1$.
This implies that the common expectation $\bE(\omega)$ of the $\omega_j$ lies in $(0,1)$,
and we claim 
that $d$ and $\bE(\omega)$ alone control our constants. We will see in Theorems \ref{t3.1} and \ref{t1.15}
that both numbers $N^E(\mu)$ and $N_u^E(\mu)$ are related to the behavior 
of the distribution function 
$F(\delta)$, and in particular its asymptotics when $\delta$ tends to $0$, which may be complicated; 
here we say that the constants in these relations depend only on $d$ and $\bE(\omega)$.

We underline  -- yet again --  that 
Theorem \ref{t1.11}
is not an asymptotic result, and multiple numerical experiments \cite{ADFJM-SIAM}  show the strength of this estimate in the intermediate regime where  $\mu$ is neither large nor small, as well as 
its
applicability to the potentials where $V$ is disordered but unbounded and 
thus, no other results for large $\mu$ are readily available. Moreover, even in the asymptotic regimes, \eqref{eq3.33-bis} offers more precision than the traditional Lifschitz tail estimates, in particular, encompassing 
faithfully the differences between individual 
choices of the disordered potentials; 
this will be discussed more thoroughly in Section \ref{S3}; 
also see~\cite{Desforges} for a detailed numerical study of the Landscape 
Law and its comparison to the available results  in the presence of disorder. 
In conclusion,  
we would like to zoom back out from the specific applications and to reiterate that the Main Theorem should be viewed as a form of the Uncertainty 
Principle  whose generality  is not inhibited by properties of the potential  or range of the energies, a ``black box" which  gives 
good
bounds on the density  of states irrespectively of the physical nature of 
the initial system.

\vskip 0.08 in

\noindent{Acknowledgements.} 
We thank Douglas Arnold and David Jerison for uncountable inspiring conversations on the subject and the joint work \cite{ADFJM-PRL, ADFJM-SIAM} which lies at the foundation of the results in this paper. 
The third author would also like to thank T. Spencer and L. Pastur for many stimulating discussions, and W. K\"onig, Z. Shen, and W. Kirsch, for sharing some references and the historical perspective. 

David is supported in part by the 
H2020 grant GHAIA 777822, and Simons Foundation grant 601941, GD.
Filoche is supported in part by Simons Foundation grant 601944, MF.
Mayboroda is supported in part by the NSF grants DMS 1344235, DMS 1839077,  and Simons Foundation grant 563916, SM.

\section{Main estimates: doubling and non-doubling scenario} 

We start with the {\it Proof of Theorem~\ref{t2.1}.} 

\noindent  {\bf Step I: the upper bound}. 
We start with the upper bound on $N(\mu)$. 
The estimate $N(\mu)\leq N$ is valid if we can find $H_N$, a codimension $N$ subspace of $H$ 
(where $H$ is the space of periodic functions in $W^{1,2}(\Omega)$),  such that 
$$
\frac{\langle Lv, v\rangle}{\|v\|_{L_2(\Omega)}^2}: = \frac{\int_{\Omega}|\nabla v|^2+V\,v^2\, dx}{\int_{\Omega}v^2\, dx} > \mu \quad \mbox{ for 
all } v\in H_N.$$
To this end, denote 
$$\F:=  \left\{Q\in \{Q\}_{\kappa\,(C_4\mu)^{-1/2}} \mbox{ such that } \min_Q \frac 1u \leq C_4 \mu\right\},$$
with $C_4$ to be defined below, and $1\leq \kappa<2$ (depending on $\mu$) 
is the smallest number such that  $R_0$ is an integer multiple of $\kappa 
\mu^{-1/2}$. Then  let $H_N$ be the space of $v\in H$ such that $\int_Q v\, dx=0$ for every $Q\in \F.$ Since the cubes $Q\in \F$ are disjoint, it is evident that $H_N$ has co-dimension $N={\rm Card }\, \F$, simply taking the bumps on $Q$'s as an orthogonal complement of $H_N$. 

We recall from \cite{ADFJM-CPDE}, Lemma~4.1,  that 
$$\int_{\Omega}|\nabla v|^2+V\,v^2\, dx\geq \int_{\Omega}\frac 1u \,v^2\, 
dx \quad \mbox{ for all } v\in H$$
and hence, 
$$2 \int_{\Omega}|\nabla v|^2+V\,v^2\, dx\geq \int_{\Omega}|\nabla v|^2+\frac 1u \,v^2\, dx \quad \mbox{ for all } v\in H.$$
Thus, it is enough to prove that 
\begin{equation}\label{eq2.7} \int_{\Omega}|\nabla v|^2+\frac 1u \,v^2\, dx >  2\mu\,\int_{\Omega}v^2\, dx   \quad \mbox{ for all } v\in H_N \sm \{0 \} . 
\end{equation}
On the part of $\Omega$ corresponding to any $Q\in \{Q\}_{\kappa\,(C_4\mu)^{-1/2}}$ such that $Q\not\in \F$ 
the bound \eqref{eq2.7} is valid provided that $C_4 >2$ because $\min_Q \frac 1u \geq C_4 \mu$ on such cubes. For $Q\in \F$, we use the Poincar\'e 
inequality to write 
$$\int_Q |\nabla v|^2\, dx\geq C_P\, C_4 \mu \int_Q |v-v_Q|^2\, dx= C_P 
\, C_4 \mu \int_Q v^2\, dx, $$
where $C_4 \mu$ comes from the size of $Q$ and we used the fact that $v_Q=\int_Q v\, dx=0$ by the definition of $H_N$. Here $C_P$ is the Poincar\'e constant and depends on the dimension only. Choosing $C_4$ so large 
that $C_P \, C_4 > 2$, we arrive at the desired estimate.

\medskip
\noindent {\bf Step II: the lower bound in  the doubling case.} 
In this direction, in order to prove that $M\leq N(\mu)$, we need to find 
$H_M$, a subspace of $H$ of dimension $M$, such that 
\begin{equation}\label{eq2.8bis}
\frac{\langle Lv, v\rangle}{\|v\|_{L_2(\Omega)}^2}: = \frac{\int_{\Omega}|\nabla v|^2+V\,v^2\, dx}{\int_{\Omega}v^2\, dx}\leq \mu \quad \mbox{ for all } v\in H_M.
\end{equation}
To this end, let 
\begin{equation} \label{deff}
\F':=  \left\{Q\in \{Q\}_{\kappa\,(C_2\mu)^{-1/2}} \mbox{ such that } \min_Q \frac 1u \leq C_2 \mu\right\},
\end{equation}
where $C_2$ will be chosen below. Let $H_M$ be the linear span of the functions $u\chi_Q$, $Q\in \F'$, picked such that $\chi_Q \in C_0^\infty(Q)$, $\chi_Q=1$ on $Q/2$, $0\leq \chi\leq 1$ on $Q$,
and $|\nabla \chi_Q|\leq 4l(Q)^{-1}$.

Since $-\Delta u\leq 1$, the Moser-Harnack inequality (\cite{HL}, Theorem 
4.14) yields 
\begin{equation}
\label{Moser-Harnack}
\sup_Q u \leq C_H \left(\frac{1}{|Q|}\int_{2Q}u^2\right)^{1/2} +C_H l(Q)^2,
\end{equation}
where $C_H$ depends on the dimension only.  In particular, using also the 
doubling condition three times, \begin{equation} \label{MH2}
\sup_{Q} u \leq C_H C_D^{3/2} \sup_{Q/4} u + C' \Bk l(Q)^2,
\end{equation}
 where $C'=C'(C_D, C_H)$ is a constant depending on $C_D, C_H$, and the 
dimension only. 

We use \eqref{eq1.2}, the definition of $\chi_Q$, 
\eqref{Moser-Harnack} for $Q/4$, and \eqref{MH2} 
\begin{multline}\label{eq2.8}
\frac{\langle L(u\chi_Q), u\chi_Q\rangle}{\|u \chi_Q\|_{L_2(\Omega)}^2} 
= \frac{\int u^2 |\nabla \chi_Q|^2 + u\chi_Q^2\, dx}{\int (u\chi_Q)^2\, 
dx}
\leq \frac{ 16\, l(Q)^{-2}\int_{Q} u^2 \,dx  + \int_{Q} u \, dx}{\int_{Q/2} u^2\, dx}\\[4pt]
\leq \frac{ 16\, l(Q)^{-2}\sup_Q u^2  +\sup_Q u}{4^{-d} \left(\frac{1}{C_H} \sup_{Q/4} u- \frac{1}{16}l(Q)^2\right)^2} 
\leq \frac{ 4^{d+2} 
\, l(Q)^{-2}\sup_Q u^2  + 4^d 
\sup_Q u}
{\left(\frac{1}{C_H^2 C_D^{3/2}} \sup_{Q} u- \bigl(\frac{1}{16}+   \frac{C'}{C_H^2 C_D^{3/2}} \bigr)l(Q)^2\right)^2} \, .
\end{multline}

We temporarily choose $\kappa$ small enough in terms of $C_D$ and $C_H$ so that  
$$\frac{1}{2 C_H^2 C_D^{3/2}} \sup_{Q} u \geq  \left(\frac{1}{16}+   \frac{C'}{C_H^2 C_D^{3/2}}  \right)l(Q)^2,$$
and then, for some constants $C'_{d,5}$, $C''_{d,5}$, $C_{d,5}$ depending 
on the dimension, $C_D$, and $C_H$ we have
\begin{equation}\label{eq2.9}
\frac{\langle L(u\chi_Q), u\chi_Q\rangle}{\|u \chi_Q\|_{L_2(\Omega)}^2}
\leq C'_{d,5} \, l(Q)^{-2} +C''_{d,5} \,\frac{1}{\sup_Q u}
\leq C_{d,5} C_2 \mu,
\end{equation}
where the last inequality comes from the definition \eqref{deff} of $\F'$.
Having fixed $\kappa$ as above, we now choose $C_2$ such that $C_{d,5} C_2=1$ and arrive at the desired estimate. To be precise, we only showed the desired inequality on the elements of the basis of $H_M$ but since the cubes $Q$ are disjoint, we immediately get it for any element of $H_M$ as well.
The only difference with what we want is that the estimate we achieved is 
in terms of the 
cardinality of a   set $\F'$ defined with an artificially small $\kappa$. 

However, if we increase the $\kappa$ to our usual fork $1\leq \kappa<2$, 
the cardinality of the resulting set $\F$ becomes even smaller, and our basis $H_M$ has less elements than expected, as desired. 

\medskip
\noindent {\bf Step III: the lower bound in the non-doubling case.} 
Our goal, once again, is to establish \eqref{eq2.8bis} for some subspace $H_M$ of dimension $M$. 
This time, we pick any $\alpha \in (0, 1/16]$ and  consider cubes of sidelength $R = \kappa\,(C^*\alpha^{d+4} \mu)^{-1/2}$.
For $Q \in \{Q\}_R$, denote by $\check Q=Q_r$ the cube concentric with $Q$ but with the smaller sidelength $r = \alpha R
= \kappa\,(C^*\alpha^{d+2} \mu)^{-1/2}$. Now take
\begin{equation} \label{otherf'}
\F':=  \left\{Q\in \{Q\}_{R} \mbox{ such that } 
\min_{\check Q} \frac 1u \leq C^*\alpha^{d+2} \mu  \, \, \mbox{ and } \min_{Q} \frac 1u \geq C^*\alpha^{d+4}\mu  \right\},
\end{equation}
and let $H_M$ be the linear span of the functions $u\chi_Q$, $Q\in \F'$, where we pick  $\chi_Q \in C_0^\infty(Q)$, $0\leq\chi_Q\leq 1$, such that 
$\chi_Q=1$ on $2 \check Q$ and $|\nabla \chi_Q|\leq CR^{-1}$. As before, we want to estimate
\begin{equation}\label{eq2.10}
\frac{\langle L(u\chi_Q), u\chi_Q\rangle}{\|u \chi_Q\|_{L_2(\Omega)}^2} 
= \frac{\int u^2 |\nabla \chi_Q|^2 + u\chi_Q^2\, dx}{\int (u\chi_Q)^2\, 
dx}
\end{equation}
(by \eqref{eq1.2}). By definition of $\F'$, $u \leq (C^*\alpha^{d+4} \mu)^{-1}$ on $Q$, so
the numerator is at most $C^2 R^{-2} \int_Q u^2 + \int_Q u 
\leq (C^*\alpha^{d+4} \mu)^{-1}|Q| \big( C^2 \kappa^{-2} + 1\big)$.
For the denominator $D$, we first apply the Moser-Harnack inequality \eqref{Moser-Harnack} 
to $\check Q$, then the definition of $\F'$, to get that
\begin{multline*}
D \geq 
\int_{2 \check Q} u^2 \geq |\check Q| \big[C_H^{-1} \sup_{\check Q} u - \ell(\check Q)^2\big]^2
= \alpha^d |Q| \big[C_H^{-1} \sup_{\check Q} u - \alpha^2 R^2\big]^2
\\
\geq \alpha^d |Q| \big[C_H^{-1} (C^*\alpha^{d+2}\mu)^{-1} - \kappa^2 \alpha^2 
(C^*\alpha^{d+4} \mu)^{-1}\big]^2
= \alpha^d |Q| (C^*\alpha^{d+2}\mu)^{-2} [C_H^{-1} - \kappa^2]^2.
\end{multline*}
We chose $\kappa^2\leq \frac{1}{2C_H}$; then the first term dominates the
second one and the expression in \eqref{eq2.10} is bounded by 
\begin{equation}\label{eq2.11}
\frac{C_{d,6} \,(C^*\alpha^{d+4} \mu)^{-1}}{C_{d,7} \alpha^d(C^*\alpha^{d+2} \mu)^{-2}}
\leq C_{d,8} C^* \mu = \mu,
\end{equation}
provided that we choose $C^\ast = C_{d,8}^{-1}$. 
Then, using the orthogonality of the $\chi_Q$, we get that
$$
N(\mu) \geq  {\rm Card} \Big\{Q\in \{Q\}_{R} 
\, ; \,   \min_{\check Q} \frac 1u \leq C^*\alpha^{d+2} \mu
\, \, \mbox{ and } \min_{Q} \frac 1u \geq C^*\alpha^{d+4}\mu  \Big\}
\geq N_1-N_2,
$$
where 
\begin{multline*}
\qquad
N_1 = {\rm Card} \Big\{Q\in \{Q\}_{R} \, ; \, \min_{\check Q} \frac 1u \leq C^*\alpha^{d+2} \mu  \Big\},
\\
N_2 = {\rm Card} \Big\{Q\in \{Q\}_{R} \, ;  \, \min_{Q} \frac 1u  \leq C^*\alpha^{d+4}\mu  \Big\}.
\end{multline*}
Notice that the cubes $\check Q = Q_r$ in this argument are smaller and 
do not cover $\Omega$,
so $N_1$ is probably not as large as 
$N'_1  = {\rm Card} \left\{R\in \{Q\}_{r} \, ; \, \min_{R} \frac 1u \leq C^*\alpha^{d+2} \mu  \right\}$. However, keeping in mind that 
we can treat $\Omega$ as a torus, we can do the estimate above for a collection of translations of 
our cubes $Q$ by a collection of at most $C \alpha^{-d}$ small vectors $e_j$, $j\in J$, so that when
we take the cubes $Q = Q_R$ as above, the smaller cubes $\check Q + e_j$, $Q \in \{Q\}_{R}$ and $j \in J$, cover $\Omega$. 
This implies that the sum of the corresponding numbers $N_1$ is at least  

$C^{-1}  N_u(C^*\alpha^{d+2} \mu)$, where $N_u$ is defined in \eqref{eq2.2} and $C$
accounts for a slight difference between $r$ and the official radius 
$\kappa (C^*\alpha^{d+2} \mu)^{-1/2}$ associated to $C^*\alpha^{d+2} \mu$.
Let us pick a nearly optimal translation $e_j$, so that
$N_1 \geq C^{-1} \alpha^{d} N_u(C^*\alpha^{d+2} \mu)$.

Similarly, $N_2 \leq C N_u(C^*\alpha^{d+4}\mu)$, and thus by the estimate 
above 
$$
N(\mu) \geq  C^{-1} \alpha^{d} N_u(C^*\alpha^{d+2} \mu) - C N_u(C^*\alpha^{d+4}\mu).
$$
This is precisely the bound \eqref{eq2.4}.
\ep 

It is important to point out that Theorem~\ref{t2.1} does not rely on the 
condition 
$V\in L^\infty(\Omega)$ and there is no dependence in constants on $\|V\|_{L^\infty(\Omega)}$ 
or on the size of the domain 
$R_0$. This is one of the main features of our estimates. 
If instead one allows our estimates to depend on $\|V\|_{L^\infty(\Omega)}$, the situation for large $\mu$ is of course rather trivial, as both the density of states and $N_u(\mu)$ roughly behave as those for the Laplacian. In particular, there exist constants $C_5, C_2, C_4$ depending on the dimension only, such that \eqref{eq2.6} is valid for all $\mu>C_5 \|V\|_{L^\infty(\Omega)}$. 
We will use an enhanced version of this statement in the next section.

\section{Anderson-type potential}
\label{S3}

We start this section with estimates on the expectation of the counting
function $N_u(\mu)$ associated to the landscape  as in \eqref{eq2.2}.

\begin{theorem}\label{t3.1} 
Let $\Omega$ and $L= -\Delta + V$ be as in Theorem \ref{t1.11}. In particular, let $\varphi\in C_0^\infty(B_{1/10}(0))$ be such that $0\leq \varphi\leq 1$,
and set
\begin{equation} \label{3.2a}
V=V_{\omega} (x)=\sum_{j\in \ZZ^d} \omega_j \, \varphi(x-j), \quad x\in \Omega,
\end{equation}
where the $\omega_j$ are i.i.d. variables taking values $\omega_j \in [0,1]$, with a 
probability distribution 
\begin{equation} \label{3.3a}
F(\delta)=\bP\{\omega \leq \delta\}, \quad 0\leq \delta\leq 1,
\end{equation}
which is not concentrated at one point, and such that 0 is the infimum of 
its support. Denote by $N^E_u(\mu)=\bE\, N_u(\mu)$ the expectation of the counting function of the minima of $1/u$, 
as defined in \eqref{eq2.2}. 
Then there exist constants $\mu^*, c_P, \gamma_1, \gamma_2$, 
depending on the dimension and the common expectation of the random variables $\omega_j$ only, 
and constants $m, \widetilde c_P, \gamma_3, \gamma_4$, depending on the dimension only,  such that 
 \begin{equation}\label{eq3.2}
\gamma_3 \,\mu^{d/2} F(\widetilde c_P \mu)^{\gamma_4 \mu^{-d/2}}\leq N^E_u(\mu) \leq \gamma_1 \,\mu^{d/2} F(c_P \mu)^{\gamma_2 \mu^{-d/2}},
\end{equation}
whenever $\mu<\mu^*$ and $R_0>(\mu m)^{-1/2}$.  
\end{theorem}

Let us put this Theorem into the context of known results for the Lifschitz tails. On the way to  our ultimate goals, we will show the following by-product of Theorem~\ref{t1.11}.

\begin{theorem}\label{t1.15} Retain the notation and assumptions of  Theorem~\ref{t1.11}.
Then there exist constants $\mu^*, m, c_P, \gamma_1, \gamma_2$, depending
on the dimension and the expectation of the random variable only, 
and constants $\widetilde c_P, \gamma_3, \gamma_4$ depending on the dimension only, 
such that 
 \begin{eqnarray}\label{eq1.12}
\gamma_3 \,\mu^{d/2} F(\widetilde c_P \mu)^{\gamma_4 \mu^{-d/2}}\leq & N_u^E(\mu) & \leq \gamma_1 \,\mu^{d/2} F(c_P \mu)^{\gamma_2 \mu^{-d/2}},\\
\label{eq1.13}
\gamma_3 \,\mu^{d/2} F(\widetilde c_P \mu)^{\gamma_4 \mu^{-d/2}}\leq & N^E(\mu) & \leq \gamma_1 \,\mu^{d/2} F(c_P \mu)^{\gamma_2 \mu^{-d/2}}
\end{eqnarray}
whenever $\mu<\mu^*$ and $R_0>(\mu m)^{-1/2}$. 
\end{theorem}

This result, and in particular  the traditionally sought-after estimate \eqref{eq1.13}, is in itself stronger than formally known asymptotics of the density of states, particularly {\it for the continuous model}, although it is fair to say that \eqref{eq1.13} would be expected by specialists 
in the subject and perhaps could even be addressed by other methods than those in the present paper. Let us explain the situation in the currently 
available literature.

The literature devoted to Lifschitz tails is extensive, particularly if one includes Poisson and other models, 
and we do not thrive here to give a comprehensive list of references or methodology -- see, e.g., \cite{Kirsch06, Konig, PasturBook} for surveys of related results. 
Here we just provide some pointers which will highlight the novelties of \eqref{eq1.13} 
(silently passing to the limit of infinite domain and removing the superscript $E$). 

The early literature, by now considered classical, and many modern textbooks treat the case when 
$F(\delta)\geq C \delta^\beta$ for some $C, \beta>0$, and provide the asymptotics
\begin{equation} \label{dwa}
\lim_{\mu\to 0} \frac{\log |\log N(\mu)|}{\log \mu}=-\frac{d}{2},
\end{equation}
see, for instance, \cite{KirschInvitation, Simon}. 
The quantity 
$$L:= \lim_{\mu\to 0} \frac{\log |\log N(\mu)|}{-\log \mu}$$
is generally known as a Lifschitz exponent, and, in addition to the results above, it is proved in \cite{PasturBook} that 
$$\lim_{\mu\to 0} \frac{\log( -\log F(\mu))}{-\log \mu}=a>0\quad  \Longrightarrow \quad L=d/2+a.$$
 
Theorems~\ref{t1.11} and \ref{t1.15} ascertain that for any non-trivial $F$ such that 
$F(\delta) > \delta$ for $\delta >0$, we can recover the Lifschitz exponent from the behavior 
of the landscape counting function 
\begin{equation}\label{Lu}L\equiv L_u, \quad\mbox{where}\quad L_u:= \lim_{\mu\to 0} \frac{\log |\log N_u(\mu)|}{-\log \mu}
\end{equation}
(assuming for simplicity that the limit exists) and 
in particular, 
\begin{equation}\label{PasturTails}L= \frac d2+\lim_{\mu\to 0} \frac{\log( -\log F(\mu))}{-\log \mu},
\end{equation}
without any a priori restrictions on $F$. This formally recovers and generalizes the results mentioned above. 
In the context of our methods, however, such statements  lose much of the 
precision exhibited in \eqref{eq3.33-bis}, \eqref{eq1.12}, \eqref{eq1.13}.  

Indeed, the problem of \eqref{dwa} is not only, or not so much, the restricted class of the potentials to which it applies, but rather the notorious imprecision of double-logarithmic asymptotics.  The underlying method of proof in  \cite{KirschInvitation, Simon} factually gives  
$$\gamma_3 \,\mu^{d/2} F(\widetilde c_P \mu)^{\gamma_4 \mu^{-d/2}}\leq  N^E(\mu)  \leq e^{-\gamma' \mu^{-d/2}}.$$
In general, the upper bound is larger than the lower bound and does not give sufficient precision to improve the double logarithm -- see the discussion and the related conjectures in \cite{KirschInvitation}. 

This is a well-known problem. The subtle difference between refined asymptotics roughly speaking asserting that $N(\mu) \sim e^{-c\mu^{-d/2}}$ and 
those with the logarithmic correction $N(\mu) \sim e^{c \mu^{-d/2}\log \mu}$ has not been overlooked in the literature. However, the refined estimates turned out to be much more challenging. At this point they are only available in $\ZZ^d$ rather than $\RR^d$ and under various additional constraints on the probability distribution -- see \cite{Konig} and \cite{M}\footnote{We are using here the review of these results from \cite{Kirsch06}. Unfortunately, the dissertation \cite{M} has never been published and so we cannot attest to the validity of the proofs or to exact statements beyond what has been quoted \cite{Kirsch06}.}. The proofs 
pass through the parabolic Anderson model -- an approach not yet developed, to the best of our knowledge, in the context of the alloy Anderson model on $\RR^d$ considered in the present paper. And, even in $\ZZ^d$, the situation has been far from well-understood. Both the conditions on the potential and the results in \cite{Konig} and \cite{M} are quite technical, so we will not provide the detailed statements.  Let us just mention that they appeal to various cases according to the behavior of the scale function 
$$S(\lambda, t)=(\lambda t)^{-1}G(\lambda t)-t^{-1} G(t), \quad \mbox{where}\quad G(t)=\log \bE(\exp (-tV(0)),$$  
(whether $S\sim C(\lambda^\rho-1)t^\rho$ with $C,\rho$ positive or negative, or $S\sim C\,\log \lambda$, or $S\sim -C(\lambda t)^{-1}\, \log t$) and draw the asymptotics in terms of of $I(\mu)=\sup_{t>0} (\mu t-G(t)).$ Such is the presentation in \cite{M}, and \cite{Konig} gives somewhat different statements, also with a dependence on the features of a certain implicitly defined scale function. The strength of these results compared 
to Theorem~\ref{t1.15} is that, at least in some cases, they provide actual asymptotics rather than the estimates from above and below and feature 
a number of cases that we did not explicitly consider, such as unbounded potentials. The weakness is that their coverage does not encompass all potentials, even among the bounded ones, and at this point is completely restricted to $\ZZ^d$.

By contrast, Theorem~\ref{t1.15} provides a simple and universal law, covering all bounded potentials at once, clearly identifying the source of the logarithmic correction, the ``Pastur tails" \eqref{PasturTails}, the exact transition from the classical to quantum regime.  Below are just a few examples of applications of \eqref{eq1.13}: \begin{enumerate}
\item $V$ is a Bernoulli potential: $\omega$ takes values $0$ or $1$ with 
probability $1/2$. Then
$$\gamma_3 \,\mu^{d/2} e^{-\gamma_4 \mu^{-d/2}}\leq  N^E(\mu)  \leq \gamma_1 \,\mu^{d/2} e^{-\gamma_2 \mu^{-d/2}}.
 $$
\item $V$ is given by a uniform distribution on $[0,1]$ or any other $F$ such that $F(\delta)$ is bounded from above and below by some positive power of $\delta$. This leads to logarithmic correctors predicted in the physics literature \cite{LN, PS}
$$\gamma_3 \,\mu^{d/2} e^{\gamma_4 \, \mu^{-d/2} \log \mu }\leq  N^E(\mu) 
 \leq \gamma_1 \,\mu^{d/2} e^{\gamma_2\, \mu^{-d/2} \log \mu }.$$

\item $V$ is given by the probability distribution with $F(\delta)\sim e^{-C\delta^{-a}}$, $a>0$. Then
$$\gamma_3 \,\mu^{d/2} e^{\gamma_4 \, \mu^{-d/2-a}}\leq  N^E(\mu)  \leq \gamma_1 \,\mu^{d/2} e^{\gamma_2\, \mu^{-d/2-a}}.$$
This is an example of \eqref{PasturTails}.
\end{enumerate}

With this, we return to the proof of Theorem~\ref{t3.1}. Our initial lemma is purely deterministic.

\begin{lemma}\label{l3.4} 
Let $\Omega$ and $L$ be as in Section~\ref{intro}, with $V$ defined as follows.
Let $\varphi\in C_0^\infty(B_{1/10}(0))$  be such that $0\leq \varphi\leq 
1$, and set 
$$
V=V_{\omega} (x)=\sum_{j\in \ZZ^d} \omega_j \, \varphi(x-j), \quad x\in \Omega,
$$
where the sequence $\omega = \{\omega_j \}_{j\in \ZZ^d}$ takes values in $[0,1]$.
For $r \in [\sqrt d, R_0]$,  where we recall that $R_0$ is the scale of $\Omega$, let us denote by $Q = Q_{r}$ the maximal cube consisting of unit cubes centered on $\ZZ^d$ 
(and with edges parallel to the axes)  which is contained in $B_{r/2}(0)$. 
Since $r \geq \sqrt d$, $Q_r$ contains at least one unit cube.

Assume that $r \in [3\sqrt d, R_0]$ is such that
\begin{equation}\label{eq3.13-bis}
{\rm Card}\, \big\{j\in Q_{r}\cap \ZZ^d:\, \omega_j\geq c_P r^{-2} \big\}
\geq \lambda\, |Q_{r}|,
\end{equation}
for some $c_P, \lambda>0$. 

If $c_P$ is large enough, depending on $\lambda$ and the dimension only, then 
there exist $\eps=\eps(\lambda, d)>0$ (small) and  $M=M(\eps, \lambda, d)>0$ (large) such that if $\xi_0 \in \ol B_{r/3}(0)$ is such that
\begin{equation}\label{eq3.5}
u(\xi_0)\geq Mr^2
\end{equation} 
then 
\begin{equation}\label{eq3.7}
u(\xi) \geq (1+\eps)\,  u(\xi_0)\quad \mbox{for some point}\quad 
\xi\in \overline{B_{\sqrt{1+\eps}\,r}(\xi_0)}. \end{equation}
\end{lemma}

Again this is a deterministic statement, for which we do not care where the $\omega_j$ are coming from and probabilistic considerations are irrelevant. That is, at this point
$V$ could be any realization, even extremely unlikely, of the construction 
described in Theorem \ref{t3.1}, even if we intend to show later that our 
assumption \eqref{eq3.13-bis} 
is quite probable in some circumstances.

Here we gave a statement for a point $\xi_0 \in \ol B_{r/3}(0)$ so that we can take
$Q_{r}$ centered at the origin, but a similar statement for any $\xi_0 \in \Omega$ would be easy to obtain, 
because we could use the translation invariance of our problem by $\ZZ^d$ 
to
apply the result to $\xi_0 - \ol\xi_0$, where $\ol\xi_0 \in \ZZ^d$ is such that 
$\xi_0 - \ol\xi_0 \in \ol B_{r/3}(0)$; we assumed $r \geq 3\sqrt d$ only to guarantee
that we can find $\ol\xi_0$. We will use this comment about other centers 
$\xi_0$ later in the proof.

\vskip 0.08 in

\bp Because of the periodic nature of $\Omega$ and $L$, we may assume that $\Omega$
is centered at the origin; we do not assume that $\xi_0 = 0$ because $\ZZ^d$ plays a special
role in the definition of $V$.
\vskip 0.08 in
\noindent {\bf Step I}.
Let $\xi_0 \in \Omega$ be given, set $B_\rho = B_\rho(\xi_0)$ 
(for computations on $u$, we like to think that $\xi_0$ is the origin)
and denote by $m(\rho) = m(\xi_0,\rho)$ the average of $u$ on the sphere centered at $\xi_0$ 
with radius $\rho$. 
That is, when $d \geq 2$ we set 
$$
m(\rho)=\fint_{\partial B_\rho} u\, d\sigma, \quad \rho>0, 
$$
where $d\sigma$ is the $(d-1)$ dimensional surface measure on $\partial B_r$, 
and when $d=1$ 
$$
m(\rho)=\frac{u(\xi_0+\rho)+u(\xi_0-\rho)}{2}, \quad \rho>0. 
$$
For brevity, we set $m(0) = u(\xi_0)$; this makes sense because $u$ is continuous on $\Omega$.
We claim that 
\begin{equation}\label{eq3.8}
m(\rho)\leq m(r)+r^2-\rho^2\quad\mbox{for}\quad 0\leq \rho<r < \dist(\xi_0,\d\Omega),  
\end{equation}
and in particular, 
\begin{equation}\label{eq3.9}
m(r)\geq m(0)-r^2. 
\end{equation}
This can be seen, for instance, by comparison with harmonic functions. Let $v$ be a solution to $-\Delta v=0$ in $B_r$ that coincides with $u$ on $\partial B_r$ and set $w(y):=v(y)+r^2-|y-\xi_0|^2$ for 
$y\in B_r$. Then $-\Delta w=2d\geq 1 \geq -\Delta u$ in $B_r$ 
(because $L u = -\Delta u + Vu = 1$ and $V \geq 0$)  and $w=v=u$ on $\partial B_r$. Hence, $w\geq u$ by the maximum principle, so that  $$m(\rho)=  \fint_{\partial B_\rho} u\, d\sigma\leq \fint_{\partial B_\rho} w\, d\sigma =\fint_{\partial B_\rho} v\, d\sigma+r^2-\rho^2=m(r)+r^2-\rho^2,$$
where we used the mean value property for harmonic functions in the last equality. The estimates \eqref{eq3.8}--\eqref{eq3.9} follow.

Furthermore, when $d \geq 2$,  the Poisson formula for a harmonic function $v$ in $B_r$ yields $$v(y)= \frac{r^2-|y-\xi_0|^2}{d\alpha_d r}\int_{\partial B_r}\frac{v(z)}{|z-y|^d}\,d\sigma_z \, ,$$
where $\alpha_d$ is the volume of a unit ball in $\RR^d$. 
Hence there exists a dimensional constant $c_1$ such that $v(y)\leq c_1 m(r)$ for all $y\in B_{2r/3}$. 
The same is of course true when $d=1$, because harmonic functions on $\RR$ are affine.
Moreover, since 
$$u(y) \leq w(y)=v(y)+r^2-|y-\xi_0|^2,$$ 
we get that
\begin{equation}\label{eq3.10}
u(y)\leq c_1 m(r)+r^2 \quad\mbox{for }\, y\in B_{2r/3}.
\end{equation}
Notice that $c_1$ can be taken equal to 1 when $y=x$, according to \eqref{eq3.9}.


\vskip 0.08 in \noindent {\bf Step II}. 
Now we want to use the size of $V$.
Integrating by parts against the Green function in a ball, we get for $d\geq 3$
\begin{multline}\label{eq3.11}
m(r)=m(0)+c_2 \int_{B_r} \Delta u(y) \left(|y-\xi_0|^{2-d}-r^{2-d}\right)\, dy\\
=m(0)+c_2 \int_{B_r} (Vu-1) \left(|y-\xi_0|^{2-d}-r^{2-d}\right)\, dy
\end{multline}
for some dimensional constant $c_2>0$ and as usual assuming 
$\ol B_r\subset \Omega$. 

Now assume that $0\leq r\leq R$ and $\ol B_R\subset \Omega$, and 
subtract \eqref{eq3.11} for $R$ from this; we get that 
\begin{multline}\label{eq3.12}
m(R)-m(r) =c_2 \int_{B_R\setminus B_r} (Vu-1)\left(|y-\xi_0|^{2-d}-R^{2-d}\right)\, dy\\
+c_2 \int_{B_r} (Vu-1) \left(r^{2-d}-R^{2-d}\right)\, dy
\end{multline}
Recall that we are interested in $\xi_0 \in \ol B_{r/3}(0)$, so that since $Q_r \subset B_{r/2}(0)$,
it is contained in $B_r = B_r(\xi_0)$.
We will only keep the contribution of $V$ on $Q_r$ (because we want to use its simpler structure), and since 
$$
\int_{B_R\setminus B_r} \left(|y-\xi_0|^{2-d}-R^{2-d}\right)\, dy + 
\int_{B_r} \left(r^{2-d}-R^{2-d}\right)\, dy 
\leq C R^d (r^{2-d}-R^{2-d}) 
\leq C (R^2-r^2)
$$
\eqref{eq3.12} yields
\begin{equation} \label{3.12b}
m(R)-m(r) \geq -c_3 (R^2-r^2)+c_2 \left(r^{2-d}-R^{2-d}\right) \int_{Q_r}Vu \, dy,
\end{equation}
In dimension $d=2$ one has 
\begin{equation}\label{eq3.11-1} m(r)=m(0)+c_2 \int_{B_r} (Vu-1) \log \frac{r}{|y-\xi_0|}\, dy
\end{equation}
in place of \eqref{eq3.11}, and since
$$
\int_{B_R\setminus B_r} \log \frac{r}{|y-\xi_0|} dy + \int_{B_r} \log \frac{R}{r} dy 
\leq C (R^2-r^2) +  C r^2 \log \frac{R}{r} \leq C (R^2-r^2),
$$
we obtain
\begin{equation}\label{eq3.11-2} m(R)-m(r) \geq  -c_3 (R^2-r^2)
+c_2 \log\frac Rr\,  \int_{Q_r}Vu \, dy
\end{equation}
in place of \eqref{3.12b}. In dimension $d=1$, \eqref{eq3.11} becomes 
\begin{equation}\label{eq3.11-3} 
m(r)= m(0) + c_2 \int_{B_r} u''(y) \,(r-|y-\xi_0|)\, dy
= m(0)+c_2 \int_{B_r} (Vu-1)\,(r-|y-\xi_0|)\, dy
\end{equation}
and hence we have 
\begin{equation}\label{eq3.11-4} m(R)-m(r) \geq  -c_3 (R^2-r^2)+c_2 \,(R-r)  \int_{Q_r}Vu \, dy
\end{equation}
in place of \eqref{3.12b}.

\vskip 0.08 in \noindent {\bf Step III}. 
Write $Q_r=\bigcup_{j\in J} R_j$,
where $R_j$ is the cube of unit sidelength centered at 
$j\in \ZZ^d$, and $J = \ZZ^d \cap Q_r$ precisely corresponds to the cubes $R_j$ that are
contained in $Q_r$. Then set
\begin{equation} \label{3.17g}
J_V:=\big\{j\in J \, : \, \omega_j\geq c_P r^{-2}\big\}.
\end{equation}

Observe that since $V(x) = \sum \omega_j \,\varphi(x-j)$, with $\varphi\in C_0^\infty(B_{1/10}(0))$,
we have that $\fint_{R_j} V = \omega_j  \fint_{R_0} \varphi$, where $R_0$ (exceptionally) 
denotes the unit cube centered at $0$. Thus
\begin{equation} \label{3.18g}
J_V:=\big\{j\in J \, : \,  \fint_{R_j} V \geq c'_P r^{-2}\big\},
\end{equation}
with $c_P' = c_P \fint_{R_0} \varphi$.

Denote by $m_r$ the average of $u$ on the ball $B_r(\xi_0)$ 
(notice the difference with $m(r)$ which is an average on the sphere) and 
let $u_j:=\inf_{R_j} u$. Now pick some $\eta>0$ (a dimensional constant 
to be chosen below) and let 
\begin{equation} \label{3.19g}
J_\eta=\{j\in J_V:\, u_j<\eta\, m_r\}.
\end{equation}

\vskip 0.08 in \noindent {\bf Step IV}. We start with the case when 
$${\rm Card}\,J_\eta \geq \frac{\lambda}{2}\, |Q_r|. $$ 

By Harnack's inequality at scale 1 (see, \cite{GT}, Theorem~8.18), 
$$
\fint_{R_j} u\, dx \leq 2^d \fint_{2R_j} u\, dx 
\leq C \left(\inf_{R_j} u +1\right). 
$$
Since, in addition, $u\geq 1$ on $\Omega$  (recall that $0 \leq V \leq 1$ 
here, 
and see \cite{ADFJM-CPDE}, Proposition~3.2),  we have 
$$\fint_{R_j} u\, dx \leq C_H' \inf_{R_j} u, $$
for some constant $C_H'$ depending on the dimension only. Therefore, 
$$\int_{R_j} u\, dx=\fint_{R_j} u\, dx\leq C_H' \eta m_r \quad \mbox{for any } \, 
j\in J_\eta.$$ 
Then 
$$
\int_{B_r\setminus \bigcup_{j\in J_\eta} R_j} u\, dx \geq |B_r| m_r - C_H' \eta m_r \,{\rm Card}\,J_\eta
$$
and 
$$ \fint_{B_r\setminus \bigcup_{j\in J_\eta} R_j} u\, dx 
\geq \frac{|B_r|  - C_H' \eta  \,{\rm Card}\,J_\eta}{|B_r|-{\rm Card}\,J_\eta} \,m_r \\
=\Big(1+ \frac{(1-C_H' \eta)  \,{\rm Card}\,J_\eta}{|B_r|-{\rm Card}\,J_\eta} \Big) m_r 
\geq (1+c_3 \lambda)\, m_r,
$$
for $\eta=(2\,C_H')^{-1}$ and a suitable dimensional constant $c_3$. 
We conclude that there exists a point $\xi\in B_r$ such that 
\begin{multline}\label{eq3.13} 
u(\xi)\geq (1+c_3 \lambda)\, m_r \geq (1+c_3\lambda) (m(0)-r^2) \\ 
\geq m(0) +c_3\lambda m(0)-m(0)\,(1+c_3\lambda)/M 
\end{multline}
where we integrated \eqref{eq3.9} for the second inequality and used 
the fact that $m(0) = u(\xi_0) \geq M r^2$ by  
\eqref{eq3.5} in the third one. If we fix 
\begin{equation}\label{eq3.14}
M\geq \frac{c_4}{\lambda}
\end{equation}
then there exists a point $\xi\in B_r$ such that 
\begin{equation}\label{eq3.15}
u(\xi) \geq (1+c_5\lambda) \,m(0),
\end{equation}
where as usual all $c_i$ depend on the dimension only. Hence, choosing 
\begin{equation}\label{eq3.16}
\eps<c_5\lambda,
\end{equation}
we arrive at \eqref{eq3.7}.

\vskip 0.08 in \noindent {\bf Step V}. Assume now that, on the contrary, 
$${\rm Card}\,J_\eta \leq \frac{\lambda}{2}\, |Q_r|. $$ 
Let $R=\sqrt{1+\eps}\, r$, $\eps<1/2$.
First assume that $d \geq 3$; then by \eqref{3.12b}, 
\begin{multline*}
m(R)-m(r) \geq -c_3 (R^2-r^2)+c_2 \left(r^{2-d}-R^{2-d}\right) \int_{Q_r}Vu \, dy\\
\geq -c_3 (R^2-r^2)+c_2 \left(r^{2-d}-R^{2-d}\right) 
\sum_{j\in J_V\setminus J_\eta}\int_{R_j}Vu \, dy.
\end{multline*}
But for such $j$, $\int_{R_j}Vu \, dy \geq u_j \int_{R_j}V \, dy \geq \eta m_r \int_{R_j}V
=\eta m_r \fint_{R_j}V \geq \eta m_r c'_P r^{-2}$ by various definitions including \eqref{3.18g}
and \eqref{3.19g}. Thus, since $R=\sqrt{1+\eps}\, r$,
\begin{multline*} 
m(R)-m(r) \geq -c_3 (R^2-r^2)
+c_2 \left(r^{2-d}-R^{2-d}\right) c_P' r^{-2} \eta m_r \,({\rm Card}\, J_V-{\rm Card}\, J_\eta)
\\
\geq -c_3 \eps r^2+c_6\eps c_P' m_r \,\lambda.
\end{multline*}
When $d=1,2$, we use \eqref{eq3.11-2} and \eqref{eq3.11-4} instead of \eqref{3.12b}, and get the same final estimate, namely
$$m(R)-m(r) \geq -c_3 \eps r^2+c_6\eps c_P' m_r \,\lambda 
$$
(possibly further adjusting $c_3$ and $c_6$ still depending on dimension only).
Using \eqref{eq3.9} and its integrated version for $m_r$, 
and then the fact that $m(0) \geq M r^2$ by \eqref{eq3.5}, 
we obtain that 
\begin{multline*}
m(R)\geq m(0)-r^2  -c_3 \eps r^2+c_6\eps c_P' \,\lambda \big(m(0)-r^2 \big)
\\
\geq m(0) \Big(1 + c_6\eps c_P'  \,\lambda \big(1-\frac 1{M}\big)  
-\frac {1+c_3\eps}{M}\Big).
\end{multline*}
Choosing $c_P$ so large that 
\begin{equation}\label{eq3.17}
c_P' \geq \frac{4}{c_6\lambda}
\end{equation} 
(recall Step III) and $M$ such that 
\begin{equation}\label{eq3.18} 
M>c_7\max\left\{1, \frac 1 \eps, \frac 1\lambda\right\}
\end{equation} 
(the third part takes care of \eqref{eq3.14})
we ensure that the second term in the parentheses above is larger than $2\eps$ 
and the third term smaller than $\eps$, so that 
$$m(R)\geq m(0) \left(1 + \eps \right) $$
and hence, \eqref{eq3.7} holds for some $\xi\in \partial B_R$, as needed for \eqref{eq3.7}.
\ep

\begin{lemma}\label{l3.19} 
Let $\Omega$ and $L = -\Delta + V$ be as in Theorems \ref{t1.11} and \ref{t3.1}.
In particular $V$ is a random potential governed by a probability measure, as in
\eqref{3.2a} and \eqref{3.3a}.
Fix $0<\lambda<1$. Then choose $c_P=c_P(\lambda,d)$ large enough, 
$\eps=\eps(\lambda, d)>0$ small enough, and  $M=M(\eps, \lambda, d)>0$ large enough, as in Lemma~\ref{l3.4}. 

Recall that $\Omega = \RR^d/R_0 \,\ZZ^d$ and, for 
$r \in [3\sqrt d, R_0]$, let $Q_r$ denote as before the maximal cube consisting of unit cubes 
centered on $\ZZ^d$  which is contained in $B_{r/2}(0)$. Then let 
\begin{equation}\label{eq3.20}
\bP_r:= \bP\big(\big\{{\rm Card}\, \{j\in Q_r \cap \ZZ^d \, :\, \omega_j\leq c_P r^{-2}\}
\geq (1-\lambda)\, |Q_r| \big\}\big).
\end{equation}
Also define a similar quantity for the whole domain, i.e.,
\begin{equation}\label{eq3.20-bis}
\bP_{\Omega} =
\bP\Big(\big\{{\rm Card}\, \{j\in \Omega \cap \ZZ^d \, :\, \omega_j\leq c_P R_0^{-2}\}
\geq  \bigl(1-\lambda\bigr)\, |\Omega| \Big\}\Big). 
\end{equation}
Finally, for $3\sqrt d \leq r < R_0$, set $r_k = (1+\eps)^{k/2} \,r$ for $0 \leq k \leq k_{max}$, 
where $k_{max}$ is the largest integer such that $r_k < R_0$. Then
\begin{equation}\label{eq3.21}
\bP \Big(\Big\{ \sup_{\xi \in \ol B_{r/3}(0)} u(\xi) \geq Mr^2 \Big\} \Big)
\leq \bP_{\Omega} + C \varepsilon^{-d} \sum_{0 \leq k \leq k_{max}}
\bP_{r_k},
\end{equation}
where $C$ depends only on the dimension.
\end{lemma}

Here we shall not even need our assumption that the probability distribution 
$F$ of \eqref{3.3a} is not concentrated at one point and $F(\delta) > 0$ for $\delta > 0$;
we will evaluate the probabilities later. 

We wrote our estimates with all the cubes $Q_\rho$, and our test ball $\ol B_{r/3}(0)$, 
all centered at $0$, but since the $\omega_j$ are i.i.d. variables and our problem is invariant 
under translations by $\ZZ^d$, the various probabilities mentioned in the 
statement would 
be the same with all the cubes (and the test ball) centered anywhere else 
on $\ZZ^d$. 
We will also use this invariance during the proof.
\smallskip

\bp The idea is to repeatedly use Lemma~\ref{l3.4} and stop when the resulting ball exceeds 
the size of $\Omega$. 

Let $r$ be given, suppose that $\sup_{\xi \in \ol B_{r/3}(0)} u(\xi) \geq 
Mr^2$; 
we pick $\xi_0  \in \ol B_{r/3}(0)$ such that $u(\xi_0) \geq Mr^2$,
and try to use Lemma~\ref{l3.4} repeatedly to find points $\xi_j$ with $u(\xi_j)$
always larger. Set (for later coherence of notation) $Q_0 = Q_r$. 
One possibility is that \eqref{eq3.13-bis} fails (with this choice of $Q_r$); 
we call this event $\Ak_0$. 
But suppose not; then Lemma~\ref{l3.4} gives a point 
$\xi_1 \in {\overline{B_{r_1}}(\xi_0)}$, with $r_1 = (1+ \varepsilon)^{1/2}$ as above,
such that $u(\xi_1) \geq (1+ \varepsilon) u(\xi_0)$, as in \eqref{eq3.7}.

Notice that $u(\xi_1) \geq M r_1^2$, so we can try to apply Lemma~\ref{l3.4} again.
This time, it could be that $\xi_1 \notin \ol B_{r_1/3}(0)$, so we choose 
$\ol\xi_1 \in \ZZ^d$
such that $\xi_1-\ol\xi_1 \in \ol B_{r_1/3}(0)$, and apply the lemma after translating by $\ol\xi_1$.
We will need to be more specific later about how we choose $\ol\xi_1$, but for the moment
let us not bother.
This means that the role of $Q_r$ is now played by $Q_1 = \ol\xi_1 + Q_{r_1}$.
One possibility is that \eqref{eq3.13-bis} fails for $Q_1$; we call this event $\Ak_1$.
But we assume not for the moment, and the lemma gives a new point  
$\xi_2 \in {\overline{B_{r_2}}(\xi_1)}$ such that $u(\xi_2) \geq (1+ \varepsilon) u(\xi_1)$, 
as in \eqref{eq3.7}. Then $u(\xi_2) \geq M r_2^2$ and we can try to apply 

Lemma~\ref{l3.4} again.

We continue as long as we do not encounter an event $\Ak_k$ where \eqref{eq3.13-bis} 
fails for $Q_k$, and then we end with a last application for $k_{max}$, which gives a
point $\xi_{k_{max}+1}$ such that $u(\xi_{k_{max}+1}) \geq M r_{k_{max}+1}^2
\geq M R_0^2$. Let $\xi_\infty \in \Omega$ be such that $u(\xi_\infty) = 
|| u ||_\infty$,
and notice that $u(\xi_\infty) \geq M R_0^2$. We now try to apply Lemma~\ref{l3.4}
one last time, to the point $\xi_\infty$, but for this it will be convenient to enlarge our domain.

Suppose for definiteness that our fundamental domain $\Omega$ 
(we abuse notation a little, and give it the same name as $\RR^d / R_0 \ZZ^d$)
is the cube of sidelength $R_0$ centered at the origin; we know that, due 
to our periodic conditions, 
other choices would be equivalent, but with this choice we were able to state and prove
Lemma~\ref{l3.4} without crossing the boundary.
Pick an odd integer $N$ larger than $4 \sqrt d$, and denote by $\wt\Omega$ the cube
centered at the origin and with sidelength $NR_0$; thus $\wt\Omega$ is composed
of $\Omega$, plus a certain number of translated copies. Extend $V$ and $u$ to be
$R_0 \ZZ^d$-periodic. Then the extension of $u$ still satisfies $L u = 1$ on $\wt\Omega$, 
and by uniqueness it is the landscape function associated to
$\wt\Omega$ and periodic boundary conditions. We apply Lemma~\ref{l3.4} with this
new, larger domain, and the radius $r = 2\sqrt d R_0$, so that the corresponding
cube $Q_r$ is precisely $\Omega$. Our choice of $N$ is large enough for this to
be possible, and also we may assume, since our problem is invariant by translations from $\ZZ^d$,
that $\xi_\infty \in \ol B_{r/3}(0)$. Our last bad event $\Ak_{k_{max}+1}$ is when
\eqref{eq3.13-bis} fails for $Q_r = \Omega$, and if this does not happen, we get
a new point $\xi \in \wt\Omega$ such that $u(\xi) \geq (1+ \varepsilon) u(\xi_\infty)$.
This is impossible, because $u(\xi_\infty) = || u ||_\infty$ and $u$ takes the same values on 
$\wt\Omega$ as on $\Omega$. 

\smallskip
At this point we proved that if the event of the left-hand side of \eqref{eq3.21}
occurs (i.e., we can find $\xi_0$ as above, with $u(\xi_0) \geq Mr^2$),
then one of the bad events $\Ak_k$ occurs. What we just need to do now is 
check that
the probability of each event $\Ak_k$ is at most the corresponding term of 
the right-hand side of \eqref{eq3.21}. In particular, we do not need to check anything about
the independence of these events, we just add their probability.

In our last case  we made sure that $Q_r = \Omega$ precisely, and so this is almost the definition 
(compare \eqref{eq3.20-bis} with \eqref{eq3.13-bis}); there is a small discrepancy, 
due to the fact that since $r = 2\sqrt d R_0$ here, we should have said 

$\omega_j \leq C_P (2\sqrt d R_0)^{-2}$ rather than $\omega_j \leq C_P R_0^{-2}$, 
but the difference only amounts to making $C_P$ a little larger, which is 
not a problem, 
and we prefer the less sharp, but simpler form in \eqref{eq3.20-bis}.

For $0 \leq k \leq k_{max}$, we need to evaluate the probability of the event $\Ak_k$,
but we have to be a little careful, because we only know that \eqref{eq3.13-bis}
fails for the translated cube $Q_k = \ol\xi_k + Q_{r_k}$, but a priori we do not 
know which cube this is. Given the position of $\xi_0 \in \ol B(r/3)(0)$, 
and the fact
that for $0 \leq m < k$, $|\xi_{m+1}-\xi_m| \leq r_{m+1}$, we see that
$|\xi_k| \leq \sum_{m=1}^{k} r_m \leq C\varepsilon^{-1} r_k$. 
We need to find $\ol\xi_k \in \ZZ^d$ such that $|\xi_k-\ol\xi_k| \leq r_k/3$, so we can
choose $\ol\xi_k$ in some set $\Xi_k$, known in advance, with less than $C \varepsilon^{-d}$
elements. Our event $\Ak_k$ can only happen if \eqref{eq3.13-bis} fails for
one of the cubes $\ol\xi + Q_{r_k}$, $\ol\xi \in \Xi_k$, and the total probability
that this happens is at most $C \varepsilon^{-d} \bP(r_k)$ (all the smaller events
associated to a single $\ol\xi \in \Xi_k$ have the same probability $\bP(r_k)$,
because our $\omega_j$ are i.i.d.). This completes the proof of Lemma \ref{l3.19}.
\ep

\begin{lemma}\label{l3.23} Let $Q$ be some cube in $\RR^d$ and 
assume that the $\omega_j$, $j\in \ZZ^d\cap Q$, are i.i.d. variables taking values $0\leq \omega_j\leq 1$, with a probability distribution 
$$F(\delta)=\bP\{\omega \leq \delta\}, \quad 0\leq \delta\leq 1,$$
 which is not trivial, i.e., not concentrated at one point, and such that 
0 is the infimum of the  support.

Fix $0<\mu<1$, $c_P^*>0$, and consider $r>0$ such that $\mu-F(c_P^* r^{-2}) >0$. 
Then such that 
\begin{multline}\label{eq3.24}
\bP\left(\left\{{\rm Card}\, \big\{j\in Q\cap \ZZ^d:\, \omega_j\leq c_P^* 
r^{-2}\big\}
\geq \mu\, {\rm Card}\, \{Q\cap \ZZ^d\}\right\} \right)
\\
\leq \big(H(\mu) F(c_P^* r^{-2})^\mu\big)^{{\rm Card}\, \{Q\cap \ZZ^d\}}
\end{multline}
with $H(\mu) = \big(\mu^\mu (1-\mu)^{1-\mu}\big)^{-1}$.
\end{lemma}
While we intend to use the Lemma for $\bP_r$ and $\bP_\Omega$ from Lemma~\ref{l3.19}, we chose to state it in full generality to emphasize explicit dependence on the constants which could be useful in other contexts.  Also, observe that 
\begin{equation} \label{mumu}
\lim_{\mu \to 1} H(\mu) = 1;
\end{equation}
we will be able to choose $\mu$ so close to $1$, depending on $\bE(\omega)$ and the dimension only,
that $H(\mu) F(c_P^* r^{-2})^\mu < F(c_P^* r^{-2})^{1/2}$, at least for $r$ sufficiently large, also  depending on $\bE(\omega)$ and the dimension 
only.

\vskip 0.08in

\bp Let $P$ denote the left-hand side of \eqref{eq3.24}, and define the 
random variables $\zeta_j$ equal to 1 when $\omega_j\leq c_P^* r^{-2}$ and 0 otherwise. 
By our assumptions the $\zeta_j$ are independent and identically distributed. Furthermore,
$$
P =  \bP\Big(\Big\{\sum_{j\in Q\cap \ZZ^d}\zeta_j 
\geq \mu\, {\rm Card}\, \{Q\cap \ZZ^d\}\Big\}\Big),
$$
hence for any $t > 0$,
\begin{equation} \label{eq3.25}
P = \bP\left(\left\{e^{t\sum_{j\in Q\cap \ZZ^d} \zeta_j} 
\geq e^{t\mu\, {\rm Card}\, \{Q\cap \ZZ^d\}}\right\} \right) 
\leq e^{-t\mu\, {\rm Card}\, \{Q\cap \ZZ^d\}} A 
\end{equation}
by Chebyshev's inequality,
and where $A$ is the expectation of the product of independent 
identically distributed variables $e^{t\zeta_j}$, hence $A = A_0^{{\rm Card}\, \{Q\cap \ZZ^d\}}$, 
where $A_0$ is the expectation of any of the $e^{t\zeta_j}$. That is,
$$
A_0 = e^t \bP(\{\omega_j\leq c_P^* r^{-2} \}) + \bP(\{\omega_j > c_P^* r^{-2} \})
= e^t F(c_P^* r^{-2}) + 1-F(c_P^* r^{-2})
$$
and, by \eqref{eq3.25}, 
$$
P \leq \exp \left( -{\rm Card}\, \{Q\cap \ZZ^d\} \big(t\mu - \log A_0\big)\right)
$$
for every $t>0$. We now want to optimize in $t$, but let us introduce notation
before we compute. Set $N = {\rm Card}\, \{Q\cap \ZZ^d\}$, $F = F(c_P^* r^{-2})$
(two constants) and, for $t>0$,
$$
f(t): = t\mu - \log A_0 = t\mu -\log (e^t F + 1-F). 
$$
Thus $P \leq e^{-N f(t)}$, and we study $f$. First, 
$f(0)=0$, and $f'(t) = \mu - \frac{e^t F}{e^t F + 1-F}$.
Thus $f'(0)=\mu-F = \mu-F(c_P^* r^{-2})>0$ by our assumptions, and hence
$f$ is increasing near $0$. In fact, $f'$ only vanishes at the point $t^*$ such that 
$$e^{t^*}=\frac{\mu}{1-\mu}\, \frac {1- F}{F}$$ 
(notice that this last value is $>1$
since $\mu > F$). Since we strongly expect $f(t)$ to be minimal for $t=t^\ast$, we decide 
to take $t = t^\ast$ in the inequality above. This yields
\begin{multline}\label{eq3.27}
P \leq e^{-N f(t^\ast)}= e^{-N t^\ast\mu + N\log (e^{t^\ast} F + 1-F))}
\\
= \exp \left(-N \mu \log\Big(\frac{\mu}{1-\mu}\, \frac {1- F}{F}\Big) 
+ N \log \big(e^{t^\ast} F + 1-F \big) \right)
\\
= \exp \left(-N \mu \log\Big(\frac{\mu}{1-\mu}\, \frac {1- F}{F}\Big) 
+ N \log \Big(\frac{1-F}{1-\mu} \Big) \right)
\\
= \exp \left(-N \log\Big(
\frac{\mu^\mu (1-\mu)^{1-\mu}}{F^{\,\mu} \,(1- F)^{1-\mu}}
\Big) \right)
= \Big( \frac{F^{\,\mu} \,(1- F)^{1-\mu}}{\mu^\mu (1-\mu)^{1-\mu}}
\Big)^N.
\end{multline}
We may drop $(1- F)^{1-\mu} \leq 1$, and now this
is the same thing as \eqref{eq3.24}; Lemma \ref{l3.23} follows.
\ep

\smallskip 
\begin{corollary}\label{c3.28} 
Let $\Omega$, $L$, and $V$ be as in Theorem \ref{t3.1}.
There exist constants $R^*, c_P,  M, \gamma_1, \gamma_2$, depending only on the dimension
and the common expectation of the random variables $\omega_j$, such that 
\begin{equation}\label{eq3.29}
\bP\left\{u(\xi_0)\geq Mr^2 \right\} \leq\gamma_1  F(c_P\, r^{-2})^{\,\gamma_2\, r^{d}}
\end{equation} 
for any $\xi_0\in \Omega$ and any $r \in (R^*, R_0]$. 
\end{corollary}

\bp
This will follow from a combination of Lemmas ~\ref{l3.19} and \ref{l3.23}. First recall our assumption that the measure associated to $F$ (call it $\nu$)
is nontrivial. Let $\bE(\omega)$ denote the expectation of our random variables; then
\begin{equation} \label{expect}
0 < \bE(\omega) < 1.
\end{equation}
where the first inequality holds because $\nu$ is not a Dirac mass
at the origin, and second one holds because the support of $\nu$
touches $0$ and is contained in $[0,1]$. 

Furthermore notice that 
$\bE(\omega) = \int_{[0,1]} \delta d\nu(\delta) =  \int_{(0,1]} \delta d\nu(\delta)
\leq 1- \nu(\{ 0 \})$ 
by Chebyshev's inequality, so $F(0) = \nu(\{ 0 \}) \leq 1- \bE(\omega) < 1$. 
Clearly, $F(c_p r^{-2})$ decays as $r$ grows. We choose a value of $F(c_p 
r^{-2})$ that we don't want to exceed, half of the way between 
$1-\bE(\omega)$ and $1$, i.e., $F_0 = \frac{2-\bE(\omega)}{2} < 1$, choose (we shall see why soon) 
$a = \frac{\bE(\omega)}{2-\bE(\omega)} \in (0,1)$,
and check now that  
\begin{equation} \label{aaa}
F(a) \leq F_0 = \frac{2 - \bE(\omega)}{2}.
\end{equation}
Indeed $\bE(\omega) = \int_{[0,1]} x d\nu(x) \leq a \nu([0,a]) + \nu((a,1]
= a F(a) + 1 - F(a)$, hence $F(a)(1-a) \leq 1-\bE(\omega)$ and since 
$1-a = \frac{2-2\bE(\omega)}{2-\bE(\omega)}$, we get \eqref{aaa}.

Now let $\mu \in (3/4,1)$ be given, to be chosen soon in terms of $F_0$, very close to $1$.
Also set $\lambda = 1-\mu$ (small), and with this $\lambda$, 
define  $c_P=c_P(\lambda,d)$ large enough, as in 
Lemma~\ref{l3.4}, and choose $\eps=\eps(\lambda, d)>0$ small enough, 
and $M = M(\eps, \lambda, d)$ large enough, again as in Lemma~\ref{l3.4}.
Those choices also work for Lemma~\ref{l3.19}, so we will be able to apply these two lemmas
with these constants.

We choose $R^\ast$ so large that $c_P (R^\ast)^{-2} \leq a$; $R^\ast$ depends 
on $\lambda$ and $\mu$, but soon we will be able to choose $\mu$ (and hence, $\lambda$), that depends
only on $\bE(\omega)$ and the dimension, 
so eventually $R^\ast$ will depend only on $\bE(\omega)$ and the dimension as well. With this choice of $R^\ast$, 
and since we shall always restrict to radii $r \geq R^\ast$, we will get that
\begin{equation} \label{aaaa}
F(c_P r^{-2}) \leq F(c_P (R^\ast)^{-2}) \leq F(a) \leq F_0 : = \frac{2 - \bE(\omega)}{2}.
\end{equation}

The whole point of Lemma \ref{l3.23} was to give a bound on the probability $\bP_r$ of 
\eqref{eq3.20}, and this bound is 
\begin{equation} \label{3.40}
\bP_r \leq \big(H(\mu) F(c_P r^{-2})^\mu\big)^{N},
\end{equation}
with $N = {\rm Card}\, \{Q\cap \ZZ^d\}$. Notice that we can take $c_P^\ast = c_P$,
and the assumption that $F(c_P r^{-2}) < \mu$ is satisfied by \eqref{aaaa} if we take
$\mu > F_0$. We also take $\mu > 3/4$, so that
$F(c_P r^{-2})^{\mu-1/2} \leq F_0^{\mu-1/2} \leq F_0^{1/4}$ and use \eqref{mumu} to finally 
choose $\mu$ so close to $1$ that $H(\mu) F_0^{1/4} < 1$. This way \eqref{3.40} implies
that $\bP_r \leq F(c_P r^{-2})^{N/2}$, which will be good enough for us.

Now let $r \geq R^\ast$ be given, and let us evaluate the probability (call it $P$) of \eqref{eq3.29}.
Notice that $P$ is smaller than the probability of having $u(\xi) \geq Mr^2$
for some point of a cube $S$ of size roughly $(10 \sqrt d)^{-1}r$, say, that contains
$\xi_0$. This probability does not depend on $S$ (by invariance), and can 
be estimated as
in Lemma~\ref{l3.19}. Thus we get that
$$
P \leq \bP_{\Omega} + C \varepsilon^{-d} \sum_{0 \leq k \leq k_{max}} \bP_{r_k},
$$
with $r_k = (1+\varepsilon)^{k/2} r$. We use (the consequence of) \eqref{3.40} to estimate 
$\bP_{r_k}$, noticing that $F(c_P r_k^{-2}) \leq F(c_P r^{-2})$ and 
each set $Q_{r_k}\cap \ZZ^d$ has at least one more point than the previous one. That is,
$N_k = {\rm Card}\{Q_{r_k}\cap \ZZ^d\}$ is at least $N+k$,
where $N = {\rm Card}\{Q_{r}\cap \ZZ^d\}$. Then
$\bP_{r_k} \leq  F(c_P r_k^{-2})^{N_k/2} \leq F(c_P r^{-2})^{(N+k)/2}\leq 
F_0^{k/2}
F(c_P r^{-2})^{N/2}$.

We have a similar estimate for $\bP_{\Omega}$ (which is of the same type as $\bP_{r_k}$,
with $r_k \sim R_0$). So we can sum the geometric series, and get the more precise estimate
\begin{equation} \label{3.41}
P \leq \gamma_1 F(c_P r^{-2})^{ {\rm Card}\{Q_r\cap \ZZ^d\}/2}
\leq \gamma_1 F(c_P r^{-2})^{\gamma_2 r^d}
\end{equation}
with constants $\gamma_1$ and $\gamma_2$ that depend on $d$ and $\bP(\omega)$
(through our choice of $F_0$, $a$, $\mu$, and then the various constants that ensue, including
$\varepsilon$). As was said earlier, we can then compute $R^\ast$, depending on these constants.
Corollary \ref{c3.28} follows.
\ep

\begin{corollary}\label{c3.31} 
Let $\Omega$, $L$, and $V$ be as in  Theorem \ref{t3.1}, in particular $V$ is a random
potential governed by i.i.d. random variables $\omega_j$. 
Then there exist constants $\mu^*, M, c_P, \gamma_3, \gamma_4$, 
depending only on the dimension and the expectation of the $\omega_j$, 
\begin{equation}\label{eq3.29aa} 
N^E_u(\mu) \leq   \gamma_3 \mu^{d/2} \, F(M c_P\,\mu)^{\,\gamma_4\, \mu^{-d/2}},
\end{equation}
whenever $\mu<\mu^*$ and $R_0>(\mu M)^{-1/2}$.

\end{corollary}
\bp Recall from \eqref{eq2.2} and the statement of Theorem \ref{t1.11} that 
$$
N^E_u(\mu)= \frac{1}{|\Omega|}\times 
\bE \left\{\mbox{the number of cubes } Q\in \{Q\}_{\kappa\,\mu^{-1/2}} \mbox{ such that } \min_Q \frac 1u \leq \mu\right\},
$$
where $1\leq \kappa<2$ (depending on $\mu$) is the smallest number such that  $R_0$ is an integer multiple of $\kappa \mu^{-1/2}$. The expectation 
of the number of cubes is less than the sum of 
expectations (by the triangle inequality), so 
$$ N^E_u(\mu) \leq \frac{1}{|\Omega|} \, \frac{|\Omega|}{(\kappa \mu^{-1/2})^d} \sup_{Q\in \{Q\}_{\kappa\,\mu^{-1/2}} } \bP \left\{ \min_Q \frac 1u \leq \mu\right\}.
$$
We want to apply our estimate in \eqref{3.41}, 
coming from Lemma~\ref{l3.19}. This one gives the probability that the infimum of $\frac 1u$
on $B_{r/3}(0)$ is at most $(Mr^2)^{-1}$,
so we should take $r$ such that $(Mr^2)^{-1}=\mu$. 
Notice that $r \leq R_0$ by our condition on $R_0$. 
We get equal probabilities for integer translations of that ball, as usual, by the translation invariance
of our setting. Now each cube $Q\in \{Q\}_{\kappa\,\mu^{-1/2}}$ can be covered by less
than $C$ integer translations of of $B_{r/3}(r/3)0)$ (taken from a fixed subgrid),
and for each one the probability that $\frac 1u \leq \mu$ somewhere on the ball is estimated as in \eqref{3.41}. Therefore
$$ N^E_u(\mu) \leq 
C (\kappa \mu^{-1/2})^{-d} \gamma_1 F(c_P r^{-2})^{\gamma_2 r^d}
\leq \gamma_3 \mu^{d/2} F(Mc_P \mu)^{\gamma_4 \mu^{-d/2}},
$$
as announced.
\ep

We now give a lower bound for $N^E_u(\mu)$.

\begin{lemma}\label{c3.30} 
Let $\Omega$, $L$, and $V$ be as in Theorem \ref{t3.1}
and in the previous lemmas.
There exist constants $m,  \widetilde c_P, \gamma_5, \gamma_6$, depending 
on the dimension only, such that 
\begin{equation}\label{eq3.31}
N^E_u(\mu) \geq   \gamma_5 \, \mu^{d/2} \, F(\widetilde c_P\,\mu)^{\,\gamma_6\, \mu^{-d/2}},
\end{equation}   
whenever $\mu\leq 1$ and $R_0>(\mu m)^{-1/2}$.

\end{lemma}
\bp 
Much as above, we start observing that 
\begin{align}
N^E_u(\mu) =\,& \frac{1}{|\Omega|}\times \bE \left\{\mbox{the number of 
cubes } Q\in \{Q\}_{\kappa\,\mu^{-1/2}} \mbox{ such that } \min_Q \frac 1u \leq \mu\right\}\\
\geq\, & \frac{1}{|\Omega|}\times \sum_{Q\in \{Q\}_{\kappa\,\mu^{-1/2}}}\bP \left\{ \min_{Q} \frac 1u \leq \mu\right\}.
\end{align}

Now we recall again from from \cite{ADFJM-CPDE}, Lemma~4.1 (or \eqref{eq1.2}),  that 
$$\int_{\Omega}|\nabla f|^2+V\,f^2\, dx\geq \int_{\Omega}\frac 1u \,f^2\, 
dx,$$
for all $f$ in the space of periodic functions in $W^{1,2}(\Omega)$, and in particular for $f\in C_0^\infty (\Omega)$. We will choose $f$ to be a standard cut-off on $4C_1Q$, $C_1\geq 1$; 
that is, $f\in C_0^\infty (4C_1Q)$, $f=1$ on $C_1Q$ and $|\nabla f|\leq 
(C_1 l(Q))^{-1}$. 
We will need that 
$4C_1Q \subset \Omega$, i.e., $\Omega$ should be large enough to accommodate this. 
This is ensured by the condition $R_0>(\mu m)^{-1/2}$  if $m$ is small enough. It follows that 
\begin{multline*}   \min_{C_1 Q} \frac 1u 
 \leq \frac{1}{|C_1Q|} \left(\int_{C_1Q} \frac 1u f^2\right)
 \leq \frac{1}{|C_1Q|} \left(\int_\Omega  |\nabla f|^2 +Vf^2\right)
 \\
 \leq \frac{1}{|C_1Q|} \left(\int_{4C_1Q}  (C_1 l(Q))^{-2} +V \right)
 \leq 4^d  \Big((C_1l(Q))^{-2}+\fint_{4C_1Q} V \Big).
 \end{multline*}
We choose $C_1$ such that $4^d C_1^{-2}\leq 1/2$; then $4^d (C_1l(Q))^{-2}
\leq l(Q)^{-2}/2 = \kappa^{-1} \mu /2 \leq \mu/2$, and now
$\min_{C_1Q} \frac 1u \leq \mu/2 + 4^d \fint_{4C_1Q} V$. Therefore
$$
\bP\left\{ \min_{C_1Q} \frac 1u\leq \mu\right\} 
\geq \bP\left\{4^d\fint_{4C_1Q}V\, dx\leq 
\mu/2\right\}
\geq \bP\left\{\max_{4C_1Q}V \leq 4^{-d}\mu/2\right\}.
$$
Note that 
$$\bP\left\{ \min_{C_1Q} \frac 1u\leq \mu\right\}\le \sum_{Q'\in\, C_1Q\bigcap  \{Q\}_{\kappa\,\mu^{-1/2}}}\bP\left\{ \min_{Q'} \frac 1u\leq \mu\right\}. $$
Therefore,  
\begin{align*}
\sum_{Q\in \{Q\}_{\kappa\,\mu^{-1/2}}}\bP\left\{ \min_{C_1Q} \frac 1u\leq 
\mu\right\} \le &
\sum_{Q\in \{Q\}_{\kappa\,\mu^{-1/2}}}
\sum_{Q'\in\, C_1Q\bigcap  \{Q\}_{\kappa\,\mu^{-1/2}}}\bP\left\{ \min_{Q'} \frac 1u\leq \mu\right\}\\
\leq & \,
C_1^d\sum_{Q\in \{Q\}_{\kappa\,\mu^{-1/2}}}
\bP\left\{ \min_{Q} \frac 1u\leq \mu\right\}.
\end{align*}

Combining all of the above and using the independence of the $\omega_j$, we conclude that
\begin{align*}
\sum_{Q\in \{Q\}_{\kappa\,\mu^{-1/2}}}
\bP\left\{ \min_{Q} \frac 1u\leq \mu\right\}
\ge & \,C_1^{-d}\sum_{Q\in \{Q\}_{\kappa\,\mu^{-1/2}}}\bP\left\{  \omega_j \leq 4^{-d}\mu/2\,\, \forall\,j\in 5C_1Q\cap \ZZ^d\right\} \\
= &\, C_1^{-d}\frac{|\Omega|}{|(\kappa \mu^{-1/2})^d|}F\left(4^{-d}\mu/2\right)^{{\rm Card}\, \{5C_1Q\,\cap \ZZ^d\}},
\end{align*}
which yields the desired conclusion. \ep 

We are now finished with the proof of  Theorem~\ref{t3.1}, which is a combination of 
Corollary~\ref{c3.31} and Lemma~\ref{c3.30}. We just renamed the four $\gamma_j$,
and also renamed $M c_P$ from Lemma~\ref{c3.30} as $c_P$, but both of these constants 
depend only on $d$ and the expectation of the $\omega_j$. 

We shall now see how Theorem~\ref{t3.1} 
provides the desired estimates on the expectation of the density of states. 

\begin{theorem}\label{t3.32} 
Let $\Omega$, $L$, and $V$ be as in Theorems \ref{t1.11} and \ref{t3.1}.
Then there exist constants $C_5, C_6 >0$, depending on the dimension
and the expectation of the random variables $\omega_j$,
only and a constant $C_4>0$, depending on the dimension only, such that 
\begin{equation}\label{eq3.33}
C_5N_u^E (C_6\, \mu) \leq N^E(\mu) \leq N_u^E (C_4 \mu),
\end{equation}
for every $\mu>0$.

In particular,  there exist constants $\mu^*, m_1, 
c_P, \gamma_1, \gamma_2$, depending on the dimension and the expectation of 
the random variable only,
and constants $\widetilde c_P, \gamma_3, \gamma_4$, depending on the dimension only, such that 
 \begin{equation}\label{eq3.34}
\gamma_3 \,\mu^{d/2} F(\widetilde c_P \mu)^{\gamma_4 \mu^{-d/2}}\leq N^E(\mu) \leq \gamma_1 \,\mu^{d/2} F(c_P \mu)^{\gamma_2 \mu^{-d/2}},
\end{equation}
whenever $\mu<\mu^*$ and $R_0>(\mu m_1)^{-1/2}$. 
\end{theorem}

Notice that Theorem \ref{t3.32} is a combination of Theorem~\ref{t1.11} and the statement \eqref{eq1.12} in Theorem~\ref{t1.15}. Since the other part of Theorem~\ref{t1.15}, \eqref{eq1.13}, 
was proved in Theorem~\ref{t3.1}, both Theorems ~\ref{t1.11} and ~\ref{t1.15} will follow 
as soon as we prove Theorem \ref{t3.32}.

\smallskip
\bp The right-hand side inequality in \eqref{eq3.33} is the right-hand side inequality in \eqref{eq2.4}, hence it has been proved in Theorem~\ref{t2.1}. The proof of the left-hand side of \eqref{eq2.4} will be split into two parts, where $\mu>\mu^\sharp$ and $\mu\le \mu^\sharp$ for some suitable $\mu^\sharp$.

For the values of $\mu>\mu^\sharp$  we are 
going to proceed as for the proof of \eqref{eq2.6} in Theorem \ref{t2.1}, 
and prove that 
for any given $\mu_0$, 
\begin{equation} \label{356}
N_u(\mu) \leq N(C' \mu) \,\, \hbox{ for all } \mu>\mu_0,
\end{equation}
where $C' = (d,\mu_0)$ depends only on $\mu_0$ and the dimension.
We will essentially use the
fact that the function $u^2$ is a doubling weight. Indeed, given that $\|V\|_{L^\infty(\Omega)} \leq 1$, the Harnack inequality  (see, \cite{GT}, Theorem~8.17 and 8.18) guarantees that 
$$\sup_{Q_{2s}} u \leq C(s) \left(\inf_{Q_s} u +s^2\right).$$
Here the constant $C(s)$  depends on $s$; specifically, the examination of the proof shows that $C(s) \leq C_0^s$ for some dimensional constant $C_0$ (see the comment right after the statement of Theorem~8.20 in \cite{GT} to this effect or simply use the Harnack inequality at scale 1 roughly 
$s$ times to treat larger $s$). 
Hence, if $s$ is bounded from above by some constant depending on $d$ and 
some 
$\mu_0>0$, we have 
$$\sup_{Q_{2s}} u \leq C(d,\mu_0) \left(\inf_{Q_s} u +s^2\right).$$
Going further, we recall that $u\geq 1$ on $\Omega$  (see \cite{ADFJM-CPDE}, Proposition~3.2), so that possibly further adjusting $C(d,\mu_0)$ we have 
$$\sup_{Q_{2s}} u \leq C(d,\mu_0) \inf_{Q_s} u,$$
again 
assuming that $s$ is bounded from above by some constant depending on $d$ 
and $\mu_0$. We now follow the argument in \eqref{eq2.8bis}--\eqref{eq2.9}, except that this time we take $C_2=1$. 
Then the sidelength of the cube under consideration is $\kappa \mu^{-1/2} 
\leq 2 (\mu_0)^{-1/2}$, 
and we will be using doubling on cubes of the size at most $16\,(\mu_0)^{-1/2}$ 
(in fact, we even use smaller $\kappa$). The argument follows the same path, 
only arriving at the bound by some constant $C'(d,\mu_0) \,\mu$ in place of 
$C_{d,5}C_2 \mu$ on the right-hand side of \eqref{eq2.9}. 
Thus \eqref{356} holds: $N_u(\mu) \leq N(C'(d,\mu_0) \mu),$ for all $\mu>\mu_0$.  We can write an upper bound on $C'(d,\mu_0)\le e^{\wt C \mu_0^{-1/2}}$ explicitly, for a suitable dimensional constant $\wt C$. Note that 
$\mu_0e^{\wt C \mu_0^{-1/2}}\to \infty$ either as $\mu_0\to 0$ or as $\mu_0\to \infty$. Therefore, we choose 
\begin{align}
\mu^\sharp= \min_{\mu_0>0}\mu_0e^{\wt C \mu_0^{-1/2}}
\end{align}
and choose $\mu_0$ to attain the minimum.
In other words, 
$$N_u(C'(d,\mu_0)^{-1}\,\mu) \leq N( \mu),$$
for all $\mu>\mu^\sharp=\mu_0\, C'(d,\mu_0)$.

Now recall the first inequality in \eqref{eq2.4} of Theorem~\ref{t2.1} and fix the constants 
$C_1, C_2, C_3$ (depending on dimension only) from this inequality. 
For the $\mu^\sharp$ given as above, we claim that for a suitable choice of $\alpha<2^{-4}$, depending on dimension and 
the expectation of the $\omega_j$, and also depending on $\mu^\sharp$,
\begin{equation}\label{eq3.35} 
C_3 N_u^E (C_2 \alpha^{d+4} \mu) 
\leq \frac 12\, C_1 
\alpha^d N_u^E (C_2 \alpha^{d+2} \mu), 
\end{equation}
whenever $\mu<\mu^{\sharp}$ and $R_0>(\mu m_1)^{-1/2}$ 
(for some $m_1 >0$, that depends on the dimension and the expectation of the $\omega_j$ only). 
As we shall see, this is basically a consequence of the fact that according to Theorem~\ref{t3.1}, $N_u^E(\mu)$ is exponentially small for small $\mu$, far beating the polynomial increase 
of $\alpha^{-d/2}$. Indeed, Theorem~\ref{t3.1} says that
\begin{equation} \label{350}
N_u^E (C_2 \alpha^{d+4} \mu) 
\leq \gamma_1 \,(C_2 \alpha^{d+4} \mu)^{d/2} 
F(c_P C_2 \alpha^{d+4} \mu)^{\gamma_2 (C_2 \alpha^{d+4} \mu)^{-d/2}},
\end{equation}
provided that $(C_2 \alpha^{d+4} \mu) < \mu^*$ and $R_0>(C_2 \alpha^{d+4} 
\mu m)^{-1/2}$.
These last conditions are ensured if we take $C_2\alpha^{d+4}\mu^\sharp \leq \mu^*$
and $m_1 \leq C_2 \alpha^{d+4} m$. 
Theorem~\ref{t3.1} also says that
\begin{equation} \label{351}
N_u^E (C_2 \alpha^{d+2} \mu) \geq 
\gamma_3 \,(C_2 \alpha^{d+2} \mu)^{d/2} 
F(\widetilde c_P C_2 \alpha^{d+2} \mu)^{\gamma_4 (C_2 \alpha^{d+2} \mu)^{-d/2}}
\end{equation}
provided that $(C_2 \alpha^{d+2} \mu) < \mu^*$ and $R_0>(C_2 \alpha^{d+2} 
\mu m)^{-1/2}$,
which will hold if we take $C_2\alpha^{d+2}\mu^\sharp \leq \mu^*$
and $m_1 \leq C_2 \alpha^{d+2} m$. 

Set $F_2 = F(\widetilde c_P C_2 \alpha^{d+2} \mu)$ and 
$F_4=F(c_P C_2 \alpha^{d+4} \mu)$; if we want to prove our claim \eqref{eq3.35},
it is enough to prove that
\begin{equation} \label{352}
C_3 \gamma_1 \,(C_2 \alpha^{d+4} \mu)^{d/2} 
F_4^{\gamma_2 (C_2 \alpha^{d+4} \mu)^{-d/2}}
\leq \frac12 C_1 
\alpha^d
\gamma_3 \,(C_2 \alpha^{d+2} \mu)^{d/2} 
F_2^{\gamma_4 (C_2 \alpha^{d+2} \mu)^{-d/2}}.
\end{equation}
Take $\alpha$ so small that $c_P C_2 \alpha^{d+4} < \widetilde c_P C_2 \alpha^{d+2}$;
thus $\alpha$ depends also on the expectation of $\omega$, through $c_P$.
Then $F_4 \leq F_2$. Also choose $\alpha$ so small that 
$\widetilde c_P C_2 \alpha^{d+2} \mu^{\sharp} < \delta_0$, with
$\delta_0 = \bE(\omega)/2$. This way, if $\nu$ denotes the probability
measure defined by $F$, 
$\bE(\omega) = \int_{[0,1]} \delta d\nu(\delta) \leq \delta_0 + \nu((\delta_0,1])
= \delta_0 + 1 - F(\delta_0)$, so $F(\delta_0) \leq 1 - \bE(\omega)/2 < 
1$.
Therefore $F_4 \leq F_2 \leq 1 - \bE(\omega)/2$ in the estimates above;
now
\begin{equation} \label{353}
\frac{F_4^{\gamma_2 (C_2 \alpha^{d+4} \mu)^{-d/2}}}
{ F_2^{\gamma_4 (C_2 \alpha^{d+2} \mu)^{-d/2}}}
\leq F_2^{a \mu^{-d/2}},
\end{equation}
with $a= \gamma_2 (C_2 \alpha^{d+4})^{-d/2}-\gamma_4 (C_2 \alpha^{d+2})^{-d/2}
\geq \frac12 \gamma_2 (C_2 \alpha^{d+4})^{-d/2}$ if $\alpha \leq (\gamma_4/\gamma_2)^{1/2}$.
Thus the right-hand side of \eqref{353} is exponentially decreasing when $\alpha$ tends to $0$.
The powers of $\mu$ in \eqref{352} are the same, and the rest is polynomial in $\alpha$;
thus \eqref{352} holds for $\alpha$ small, and \eqref{eq3.35} follows.

Now we average \eqref{eq2.4} and use \eqref{eq3.35}; we get that
\begin{equation}\label{eq355}
\frac{C_1}{2}  
\alpha^d
N^E_u (C_2 \alpha^{d+2} \mu)
\leq N^E(\mu) \leq N^E_u (C_4 \mu),
\end{equation}
which is the same as \eqref{eq3.33}
(recall that we are allowed to let $C_5$ and $C_6$ depend on $\alpha$,
which is now chosen depending on $\bP(\omega)$ and $d$), 
except that we have to assume 
that $\mu<\mu^{\sharp}$ and $R_0 > (m_1 \mu)^{-1/2}$, and 
\begin{align}
\alpha=\min\left\{\, \left(\frac{\mu^\ast}{C_2\mu^\sharp}\right)^{\frac{1}{d+2}}, \, 
\left(\frac{\delta_0}{\wt c_PC_2\mu^\sharp}\right)^{\frac{1}{d+2}}, \,
\left(\frac{\gamma_4}{\gamma_2}\right)^{\frac{1}{2}}
  \right\}
\end{align}

Taken along with Theorem~\ref{t3.1}, this also automatically gives \eqref{eq3.34}. As usual, we silently redefine the constants, 
still depending on the same parameters.

 \ep

\Addresses


\begin{thebibliography}{999} \normalsize

\bibitem[ADFJM1]{ADFJM-CPDE} D.\,Arnold, G.\,David, M\,Filoche, D.\,Jerison, and S.\, Mayboroda, {\it Localization of eigenfunctions via an effective potential.} Comm. PDE, \textbf{44}:11, 1186--1216, 2019. DOI:\href{https://doi.org/10.1080/03605302.2019.1626420}{10.1080/03605302.2019.1626420}

\bibitem[ADFJM2]{ADFJM-SIAM} D.\,Arnold, G.\,David, M\,Filoche, D.\,Jerison, and S.\,Mayboroda,  {\it Computing spectra without solving eigenvalue problems.} SIAM J. Sci. Comput., 41(1), B69–B92, 2019. DOI:\href{https://doi.org/10.1137/17M1156721}{10.1137/17M1156721}

\bibitem[ADFJM3]{ADFJM-PRL} D.\,Arnold, G.\,David, D.\,Jerison, S.\,Mayboroda, and M.\,Filoche,  {\it Effective confining potential of quantum states in disordered media.} Phys. Rev. Lett., \textbf{116}, 056602, 2016. DOI: \href{https://doi.org/10.1103/PhysRevLett.116.056602}{10.1103/PhysRevLett.116.056602}

\bibitem[CKOOS]{Schulz-PRB} D.\,Chaudhuri, J.C.\,Kelleher, M.R.\,O'Brien, E.P.\,O'Reilly, and S.\,Schulz, {\it Electronic structure of semiconductor nanostructures: A modified localization landscape theory.} Phys. Rev. B \textbf{101}, 035430, 2020. \href{https://doi.org/10.1103/PhysRevB.101.035430}{10.1103/PhysRevB.101.035430}

\bibitem[CT]{ComtetPRL} A.\,Comtet and C.\,Texier, {\it Comment on ``Effective Confining Potential of Quantum States in Disordered Media''.} Phys. Rev. Lett. \textbf{124}, 219701, 2020. DOI:\href{https://doi.org/10.1103/PhysRevLett.124.219701}{10.1103/PhysRevLett.124.219701}

\bibitem[DM+]{Desforges} P.\,Desforges, S.\,Mayboroda, S.\,Zhang, G.\,David, D.\,Arnold, W.\,Wang, and M.\,Filoche, {\it Sharp estimates for the integrated density of states in Anderson tight-binding models.} \href{http://arxiv.org/abs/2010.09287}{arxiv:2010.09287} [math-ph], 2020.

\bibitem[F]{F} C.\, Fefferman, {\it The uncertainty principle.} Bull. Amer. Math. Soc. (N.S.) \textbf{9}, 129--206, 1983. DOI: \href{https://doi.org/10.1090/S0273-0979-1983-15154-6}{10.1090/S0273-0979-1983-15154-6}

\bibitem[FADJM]{ReplyPRL} M.\,Filoche, D.\,Arnold, G.\,David, D.\,Jerison, and S.\,Mayboroda, {\it Filoche et al. Reply:}, Phys. Rev. Lett., \textbf{124}, 219702, 2020. DOI: \href{https://doi.org/10.1103/PhysRevLett.124.219702}{10.1103/PhysRevLett.124.219702,}

\bibitem[FM2]{FM-PNAS} M.\,Filoche and S.\,Mayboroda. {\it Universal mechanism for Anderson and weak localization.} Proc. Natl. Acad. Sci. USA \textbf{109}(37):14761-14766, 2012. DOI: \href{https://doi.org/10.1073/pnas.1120432109}{10.1073/pnas.1120432109} 

\bibitem[FP+]{LED1} M.\,Filoche, M.\,Piccardo,  
Y.R. Wu, C.-K. Li, C. Weisbuch, and S. Mayboroda. {\it Localization landscape theory of disorder in semiconductors I: Theory and modeling.} Phys. Rev. B \textbf{95}, 144204, 2017. DOI: \href{https://doi.org/10.1103/PhysRevB.95.144204}{10.1103/PhysRevB.95.144204}

\bibitem[GT]{GT} D. Gilbarg and N. S. Trudinger, {\it Elliptic Partial Differential
Equations of Second Order.} 2nd Edition, Springer Verlag, Berlin Heidelberg 1983.

\bibitem[HL]{HL} Q.\,Han, F.\,Lin, {\it Elliptic partial differential equations.} Second edition. Courant Lecture Notes in Mathematics, 1. Courant Institute of Mathematical Sciences, New York; American Mathematical Society, Providence, RI, 2011. 

\bibitem[K]{KirschInvitation} W.\, Kirsch, {\it An invitation to random Schr\"odinger operators.} With an appendix by Fr\'ed\'eric Klopp. Panor. Synth\`eses, 25, Random Schr\"odinger operators, 1--119, Soc. Math. France, Paris, 2008.

\bibitem[KM]{Kirsch06}  W.\,Kirsch, B.\,Metzger, {\it The integrated density of states for random Schr\"odinger operators.}  Spectral theory and mathematical physics: a Festschrift in honor of Barry Simon's 60th birthday, 649--696, Proc. Sympos. Pure Math., 76, Part 2, Amer. Math. Soc., Providence, RI, 2007. DOI: \href{https://doi.org/10.1090/pspum/076.2}{10.1090/pspum/076.2 }

\bibitem[Ko]{Konig} W.\,K\"onig, {\it The parabolic Anderson model. Random walk in random potential. } Pathways in Mathematics. Birkh\"auser/Springer, [Cham], 2016. 

\bibitem[Ku]{Kurata} K.\,Kurata, {\it On doubling properties for non-negative weak solutions of elliptic and parabolic PDE. } Israel J. Math. \textbf{115}, 285--302, 2000. DOI: \href{https://doi.org/10.1007/BF02810591}{10.1007/BF02810591}

\bibitem[LN]{LN} J.M.\,Luck , Th.M.\,Nieuwenhuizen, {\it Lifshitz tails and long-time decay in random systems with arbitrary disorder.} J. Statist. Phys. \textbf{52}, 1--22, 1988. DOI: \href{https://doi.org/10.1007/BF01016401}{10.1007/BF01016401}

\bibitem[LP+]{Dirac} G.\,Lemut, M. J.\,Pacholski, O.\,Ovdat, A.\,Grabsch, J.\,Tworzyd\l{}o, and C. W. J.\,Beenakker, {\it Localization landscape for Dirac fermions.} Phys. Rev. B \textbf{101}, 081405(R), 2020. DOI: \href{https://doi.org/10.1103/PhysRevB.101.081405}{10.1103/PhysRevB.101.081405}

\bibitem[M]{M} B.\,Metzger, {\it Asymptotische Eigenschaften im Wechselspiel von Diffusion und Wellenausbreitung in zuf\"alligen Medien,} Ph.D. thesis, TU Chemnitz (2005).

\bibitem[PF]{PasturBook} L.\,Pastur, A.\,Figotin,  {\it Spectra of random and almost-periodic operators. } Grundlehren der Mathematischen Wissenschaften [Fundamental Principles of Mathematical Sciences], 297. Springer-Verlag, Berlin, 1992.

\bibitem[PL+]{LED2} M.\,Piccardo, C.-K.\,Li, Y.-R.\,Wu, J.\,Speck, B.\,Bonef, R.\,Farrell, M.\,Filoche, L.\,Martinelli, J.\,Peretti, and C.\,Weisbuch. {\it Localization landscape theory of disorder in semiconductors. II. Urbach tails of disordered quantum well layers.} Phys. Rev. B \textbf{95}, 144205, 2017. DOI: \href{https://doi.org/10.1103/PhysRevB.95.144205}{10.1103/PhysRevB.95.144205}

\bibitem[PS]{PS} A.\,Polti and T.\,Schneider, {\it Corrections to the Lifshitz tail and the long-time behaviour of the trapping problem.} EPL (Europhysics Letters) \textbf{5}:715, 1988. DOI: \href{http://doi.org/10.1209/0295-5075/5/8/009}{10.1209/0295-5075/5/8/009}

\bibitem[S1]{Shen-TAMS} Z.\,Shen, {\it Eigenvalue asymptotics and exponential decay of eigenfunctions for Schr\"odinger operators with magnetic fields.} Trans. Amer. Math. Soc. 348 (1996), no. 11, 4465--4488.

\bibitem[S2]{Shen-Duke} Z.\,Shen, {\it On bounds of $N(\lambda)$ for a magnetic Schr\"odinger operator.} Duke Math. J. \textbf{94}(3): 479--507, 1998. DOI: \href{https://doi.org/10.1215/S0012-7094-98-09420-0}{10.1215/S0012-7094-98-09420-0}

\bibitem[S79]{S79}  B.\,Simon, {\it Functional integration and quantum physics.} Pure and Applied Mathematics, 86. Academic Press, Inc. [Harcourt Brace Jovanovich, Publishers], New York-London, 1979.

\bibitem[S]{Simon} B.\,Simon, {\it Lifschitz tails for the Anderson model.} J. Statist. Phys. \textbf{38}, 65--76, 1985. DOI: \href{https://doi.org/10.1007/BF01017848}{10.1007/BF01017848}

\bibitem[TM+]{JApplPhys} T.\,Y.\,Tsai, K.\,Michalczewski, P.\,Martyniuk, C.-H.\,Wu, and Y.-R.\,Wu, {\it Application of localization landscape theory and the k $\cdot$ p model for direct modeling of carrier 
transport in a type II superlattice InAs/InAsSb photoconductor system.} J. Appl. Phys. \textbf{127}, 033104, 2020. DOI:\href{https://doi.org/10.1063/1.5131470}{10.1063/1.5131470}

\end{thebibliography}
\end{document}